\documentclass{aims}
\usepackage{amsmath}
  \usepackage{paralist}
  \usepackage{graphics} 
  \usepackage{epsfig} 
 \usepackage[colorlinks=true]{hyperref}

\usepackage{subfigure} 
\usepackage{cite}
\usepackage{graphics,color}
\usepackage{graphicx} 
\usepackage{amssymb}  

\definecolor{litegray}{rgb}{0.88,0.88,0.88}
\definecolor{daresay}{rgb}{0.55,0.55,0.55}

\hypersetup{urlcolor=blue, citecolor=red}

\def\currentvolume{X}
 \def\currentissue{X}
  \def\currentyear{200X}
   \def\currentmonth{XX}
    \def\ppages{X--XX}

\newtheorem{theorem}{Theorem}[section]

\newtheorem{lemma}[theorem]{Lemma}
\newtheorem{proposition}[theorem]{Proposition}

\theoremstyle{definition}
\newtheorem{definition}[theorem]{Definition}
\newtheorem{remark}[theorem]{Remark}

\newcommand{\R}{{\mathbb R}}
\newcommand{\rref}[1]{\mbox{(\ref{#1})}}
\newcommand{\vv}{\vspace{1mm}}
\newcommand{\vvv}{\vspace{3mm}}
\newcommand{\bd}{\vspace{0mm}\begin{displaymath}}
\newcommand{\ed}{\vspace{0mm}\end{displaymath}}
\newcommand{\be}{\vspace{0mm}\begin{equation}}
\newcommand{\ee}{\vspace{0mm}\end{equation}}
\newcommand{\bda}{\vspace{0mm}\begin{eqnarray*}}
\newcommand{\eda}{\vspace{0mm}\end{eqnarray*}}
\newcommand{\bea}{\vspace{0mm}\begin{eqnarray}}
\newcommand{\eea}{\vspace{0mm}\end{eqnarray}}
\newcommand{\beqn}{\begin{eqnarray}}
\newcommand{\eeqn}{\end{eqnarray}}

\title[Analysis of a conservation law]
      {Analysis of a conservation law modeling\\
        a highly re-entrant manufacturing system}

\author[Jean-Michel Coron, Matthias Kawski, and Zhiqiang Wang]{}


 \email{coron@ann.jussieu.fr}
 \email{kawski@asu.edu}
 \email{wzq@fudan.edu.cn}

\thanks{The first author was partially supported by the ``Agence Nationale de
la Recherche'' (ANR), Project C-QUID, number BLAN-3-139579.
The second author was partially supported by the National Science Foundation
through DMS 05-09030, the University Pierre and Marie Curie-Paris VI
and the Foundation Sciences Math\'{e}matiques de Paris.
The third author was partially supported by the Natural Science Foundation of China
grant 10701028 and the Foundation Sciences Math\'{e}matiques de Paris.}

\begin{document}
\maketitle

\centerline{\scshape Jean-Michel Coron }
\medskip
{\footnotesize
 \centerline{Institut universitaire de France and Universit\'{e} Pierre et Marie Curie-Paris VI}
   \centerline{UMR 7598, Laboratoire Jacques-Louis Lions, 4, place Jussieu}
   \centerline{Paris, F-75005 France}
} 


\medskip

\centerline{\scshape Matthias Kawski}
\medskip
{\footnotesize
 \centerline{ Arizona State University}
   \centerline{Tempe, Arizona 85287-1804, USA}
}
\medskip

\centerline{\scshape Zhiqiang Wang}
\medskip
{\footnotesize
 \centerline{Fudan University and Universit\'{e} Pierre et Marie Curie-Paris VI}
   \centerline{UMR 7598, Laboratoire Jacques-Louis Lions, 4, place Jussieu}
   \centerline{Paris, F-75005 France}
}

\bigskip

 \centerline{(Communicated by the associate editor name)}


\begin{abstract}
This article studies a hyperbolic conservation law that models
a highly re-entrant manufacturing system as encountered  in
semi-conductor production.
Characteristic features are the nonlocal character of the velocity
and that the influx and outflux constitute the control and output
signal, respectively.
We prove the existence and uniqueness of solutions for $L^1$-data,
and study their regularity properties.
We also prove the existence of optimal controls
that minimizes in the $L^2$-sense the mismatch between the actual
and a desired output signal.
Finally, the time-optimal control for a step between equilibrium states is
identified and proven to be optimal.
\end{abstract}

\section{Introduction and prior work}

This article studies optimal control problems governed by the
scalar hyperbolic conservation law
\be
\label{eq0}
\partial_t\rho(t,x)+\partial_x
\left(\lambda(W(t))\,\rho(t,x)\right)=0
\;\;\mbox{ where }\;\;
W(t)=\int_0^1\!\!\rho(t,x)\,dx,
\ee
on a rectangular domain
$[0,T] \times [0,1]$ or the semi-infinite strip
$[0,\infty) \times [0,1]$. We assume that $\lambda(\cdot)\in
C^1([0,+\infty); (0,+\infty))$ in the whole paper.

This work is motivated by problems arising in the control of
semiconductor manufacturing systems which are characterized
by their highly re-entrant character, see below for more
details. In the manufacturing system the natural control input is the
influx, which suggests the boundary conditions
\be
\label{eq0IC}
\rho(0,x)=\rho_0(x),
\mbox{ for}\;\; 0\leq x \leq 1,\;
\mbox{ and} \;
\rho(t,0)\lambda(W(t)) = u(t),
\mbox{ for}\;\; t\geq 0.
\ee
Various different choices of the space of admissible controls are
of both practical and mathematical interest, each leading to
distinct mathematical problems.
Motivated by this application from manufacturing systems, natural
control objectives are to minimize the {\em error signal} that is
the difference between a given demand forecast $y_d$ and the actual
out-flux $y(t)=\lambda(W(t))\rho(t,1)$.
An alternative to this problem modeling a {\em perishable demand},
is the similar problem that permits {\em backlogs}. In that case,
the objective is to minimize in a suitable sense the size of the
different error signal
\be
\label{opt-backlog}
\beta (t)=\int_0^t y_d(s)\, ds
         -\int_0^t \lambda (W(s)) \rho(s,1) \, ds,
\ee
while keeping the state $\rho(\cdot,x)$ bounded.
This article only considers the problem of perishable demand and
the minimization in the $L^2$-sense.
\vv

Partial differential equations models for such manufacturing systems
are motivated by the very high volume (number of parts manufactured
per unit time) and the very large number of consecutive production
steps which typically number in the many hundreds.
They are popular due to their superior analytic properties and the
availability of efficient numerical tools for simulation.
For more detailed discussions see e.g.
\cite{
DA-CR-degond,
DA-CR-herty5,
DA-CR-hier,
DA-CR-thermalized,
herty-network,
herty-supp,
lamarca08}.
In many aspects these models are very similar to those of traffic
flows, compare e.g. \cite{piccoli-cars}.
\vvv

The study of hyperbolic conservation laws,
and especially of control systems governed by such laws,
have a rich history.
A modern introduction to the subject is the text \cite{bressanbook}.
From a mathematical perspective, the choice of spaces in which to
consider the conservations laws (and their data) provides for
distinct levels of challenges.
Fundamental are question of wellposedness,
regularity properties of solutions,
controllability,
existence, uniqueness and regularity of optimal controls.
Existence of solutions, regularity and well-posedness of nonlinear
conservation laws have been widely studied under diverse sets of
hypotheses, commonly in the context of vector values systems of
conservation laws, see e.g.
\cite{fabio-exist,
piccoli2001,
bress-wellposed}.
Further results on uniqueness may be found in \cite{bress97},
while \cite{shift-diff} introduced an a distinct notion of
differentiability of the solution of hyperbolic systems.
For the controllability of linear hyperbolic systems,
see, in particular, the important survey \cite{1978-Russell}.
The attainable sets of nonlinear conservation laws are studied
in \cite{fabio-attset,coron09,Horsin,libook,Li-Rao},
while \cite{coronbook} provides a comprehensive survey of
controllability that also includes nonlinear conservation laws.
\vvv

This article is, in particular, motivated by the recent work \cite{lamarca08}
which, among others, considered the optimal control problem of minimizing
$\|y-y_d\|_{L^2(0,T)}$ (the $L^2$ norm of the difference between a demand forecast
and the actual outflux).
That work derived necessary conditions and used these to numerically
compute optimal controls corresponding to piecewise constant desired
outputs $y_d$.

The organization of the following sections is as follows:
First we rigorously prove the existence of weak solutions of the
Cauchy problem for the conservation law \rref{eq0} for the case when the initial
data and boundary condition \rref{eq0IC} lie in $L^1(0,1)$ and $L^1(0,T)$, respectively.
Next we establish the existence and uniqueness of solutions for
the optimal control problem of minimizing the $L^2$-norm of the
difference between any desired $L^2$-demand forecast $y_d$ and
actual outflux $y(t)=\lambda(W(t))\cdot \rho(t,1)$.
Finally, in the classical special case where
\be
\label{usualspeed}
\lambda(W)={1\over 1+W},
\ee
we prove that the natural candidate control for transferring
the system from one equilibrium state to another one is indeed
time-optimal.

\vv

While preparing the final version of this article, the authors
received a copy of the related manuscript \cite{herty09} which is
also motivated in part by
\cite{DA-CR-degond,DA-CR-herty5,lamarca08}
and which addresses wellposedness for systems of hyperbolic
conservation laws with a nonlocal speed on all of $\R^n$.
It also includes a study of the solutions with respect
to the initial datum and  a necessary condition for the
optimality of  integral functionals.
There are substantial differences between \cite{herty09} and our paper,
especially the treatment of the boundary conditions and the method of proof.


\section{Existence, uniqueness, and regularity of solutions in $L^1$}

\subsection{Technical preliminaries and notation}

For any $\lambda \in C^1([0,+\infty); (0,+\infty))$
define the functions $\widetilde\lambda, \overline\lambda \in
C^0([0,\infty);(0,\infty))$ and $d\in C^0([0,\infty);[0,\infty))$
with respect to $\lambda$ as
  \be \label{lambda-bound}
  \widetilde\lambda(M):=\inf_{0\leq W\leq M} \lambda(W), \
  \overline\lambda(M):=\sup_{0\leq W\leq M} \lambda(W), \
    d(M) :=\sup_{0\leq W\leq M} |\lambda'(W)|.
  \ee
 For convenience we extend $\lambda$ to all of $\mathbb{R}$ in such a way that this extension,
 still denoted $\lambda$, is in $C^1(\mathbb{R};(0,+\infty))$.

\subsection{Weak solutions of the Cauchy problem}

First we recall, from \cite[Section 2.1]{coronbook}, the usual definition of a
weak solution to the Cauchy problem \rref{eq0} and \rref{eq0IC}.

\begin{definition}
\label{weaksol}
Let $T>0$, $\rho_0\in L^1(0,1)$ and $u\in L^1(0,T)$ be given. A weak
solution of the Cauchy problem \rref{eq0} and \rref{eq0IC} is a
function $\rho\in C^0([0,T];L^1(0,1))$ such that, for every
$\tau\in[0,T]$ and every $\varphi\in C^1([0,\tau]\times[0,1])$ such
that
 \be \label{varphi}
 \varphi(\tau,x)= 0,\forall x\in[0,1]\quad \text{and}\quad
 \varphi(t,1)= 0,\forall t\in[0,\tau],
 \ee
one has
 \bea
 \label{DefSol}
 \int_0^{\tau} \int_0^1 \rho(t,x)(\varphi_t(t,x)
 +\lambda(W(t))\varphi_x(t,x)) dx dt
  \nonumber \\
  +\int_0^{\tau} u(t)\varphi(t,0)dt
  +\int_0^1 \rho_0(x)\varphi(0,x)dx&=&0.
 \eea
\end{definition}
\vv

One has the following lemma, which will be useful to prove
a uniqueness result for the Cauchy problem \rref{eq0} and \rref{eq0IC}.

\begin{lemma}\label{lem-test-function}
If $\rho\in C^0([0,T];L^1(0,1))$ is a weak solution to the Cauchy
problem \rref{eq0} and \rref{eq0IC}, then for every $\tau\in[0,T]$
and every $\varphi\in C^1([0,\tau]\times[0,1])$ such that
 \be\label{varphi-2}
  \varphi(t,1)= 0,\forall t\in[0,\tau],
 \ee
one has
 \bea \label{test}
 \int_0^{\tau} \int_0^1 \rho(t,x)(\varphi_t(t,x)
   +\lambda( W(t))\varphi_x(t,x))dxdt +\int_0^{\tau}u(t)\varphi(t,0)dt
  \nonumber\\
 -\int_0^1 \rho(\tau,x)\varphi(\tau,x)dx
    +\int_0^1\rho_0(x)\varphi(0,x)dx &=&0.
 \eea
\end{lemma}

\begin{proof}
The case $\tau=0$ is trivial. For every $\tau\in (0,T]$ and
$\varepsilon \in (0,\tau)$, let $\eta_{\varepsilon} \in
C^1([0,\tau])$ be such that
 \be \label{eta-epsilon}
  \eta_{\varepsilon}(\tau)=0
  \quad \text{and}\quad  \eta_{\varepsilon}(t)=1,\ \forall t\in[0,\tau-\varepsilon]
  \quad \text{and}\quad  \eta_{\varepsilon}'(t) \leq 0,\ \forall t\in  [0,\tau].
 \ee
It is easy to prove that, for every $h\in C^0([0,\tau])$,
 \be \label{eta-epsilon-limit}
 \lim_{\varepsilon\rightarrow 0} \int_{\tau-\varepsilon}^{\tau}
 \eta_{\varepsilon}'(t)h(t) dt=-h(\tau).
 \ee
Then, for every  $\varphi\in C^1([0,\tau]\times[0,1])$ satisfying
\rref{varphi-2}, let
$\varphi_{\varepsilon}(t,x):=\eta_{\varepsilon}(t) \varphi(t,x)$.
This obviously verifies
 \be \label{varphi-epsilon}
 \varphi_{\varepsilon}(\tau,x)= 0,\forall x\in[0,1]\quad \text{and}\quad
 \varphi_{\varepsilon}(t,1)= 0,\forall t\in[0,\tau].
 \ee

Since $\rho\in C^0([0,T];L^1(0,1))$ is a weak solution to the Cauchy
problem \rref{eq0} and \rref{eq0IC}, we have
 \bea \label{test-2}
 \lefteqn{\int_0^{\tau} \int_0^1 \rho(t,x)((\varphi_{\varepsilon})_t(t,x)
   +\lambda( W(t))(\varphi_{\varepsilon})_x(t,x))dxdt  }
  \nonumber\\
 &&  +\int_0^{\tau}u(t)(\varphi_{\varepsilon})(t,0)dt
     +\int_0^1\rho_0(x)(\varphi_{\varepsilon})(0,x)dx =0.
 \eea
Using the definition of $\varphi_{\varepsilon}$, \rref{eta-epsilon} and \rref{test-2}, one has
 \bea\label{test-3}
 \lefteqn{ \int_0^{\tau} \int_0^1 \rho(t,x)(\varphi_t(t,x)
   +\lambda( W(t))\varphi_x(t,x))dxdt}
   \nonumber\\
 && +\int_0^{\tau}u(t)\varphi(t,0)dt
    +\int_0^1\rho_0(x)\varphi(0,x)dx
  \nonumber\\
 &=& \int_{\tau-\varepsilon}^{\tau} \int_0^1
   (1-\eta_{\varepsilon}(t)) \rho(t,x)(\varphi_t(t,x) +\lambda( W(t))\varphi_x(t,x))dxdt
  \nonumber\\
 &&   + \int_{\tau-\varepsilon}^{\tau} (1-\eta_{\varepsilon}(t)) u(t)\varphi(t,0) dt
     -\int_{\tau-\varepsilon}^{\tau} \int_0^1 \eta_{\varepsilon}'(t)\rho(t,x)\varphi(t,x) dx dt.
 \eea
Observing that $\rho\in C^0([0,T];L^1(0,1))$, $\lambda\in
C^1(\mathbb{R}; (0,\infty))$ and $\varphi\in C^1([0,\tau] \times
[0,1])$, we point out that the functions $W(\cdot) =\int_0^1
\rho(\cdot,x) dx $, $\int_0^1 \rho(\cdot,x)\varphi(\cdot,x) dx $ and
$\lambda(W(\cdot))$ are all in  $C^0([0,T])$.

We can estimate the first two terms on the right hand side of
\rref{test-3} as
 \be \label{test-3-1}
  \Big| \int_{\tau-\varepsilon}^{\tau} \int_0^1
   (1-\eta_{\varepsilon}(t)) \rho(t,x)(\varphi_t(t,x)
   +\lambda( W(t))\varphi_x(t,x))dxdt \Big|
   \leq K\varepsilon,
 \ee
and  \be \label{test-3-2}
  \Big| \int_{\tau-\varepsilon}^{\tau} (1-\eta_{\varepsilon}(t))
   u(t)\varphi(t,0) dt \Big|  \leq   K\int_{\tau-\varepsilon}^{\tau} u(t)
   dt,
 \ee
where $K$ is a constant independent of $\varepsilon$. While for the
last term on the right hand side of \rref{test-3}, we get from
\rref{eta-epsilon-limit} that
 \bea \label{test-3-3}
 \lefteqn{\int_{\tau-\varepsilon}^{\tau} \int_0^1 \eta_{\varepsilon}'(t)\rho(t,x)\varphi(t,x) dx dt
 =\int_{\tau-\varepsilon}^{\tau}\eta_{\varepsilon}'(t)\Big( \int_0^1 \rho(t,x)\varphi(t,x) dx\Big) dt
 }
 \nonumber\\
  && \longrightarrow -\int_0^1 \rho(\tau,x)\varphi(\tau,x) dx \quad
  \text{as}  \quad   \varepsilon \rightarrow 0.
  \rule{50mm}{0mm}
 \eea
In view of \rref{test-3-1}-\rref{test-3-3}, letting
$\varepsilon \rightarrow 0$ in \rref{test-3} one gets
\rref{test}.
\end{proof}

\begin{theorem}\label{thm-l1}
If $\rho_0\in L^1(0,1)$ and $u\in L^1(0,T)$ are nonnegative almost
everywhere, then the Cauchy problem \rref{eq0} and \rref{eq0IC}
admits a unique weak solution $\rho\in C^0([0,T];L^1(0,1))$, which
is also nonnegative almost everywhere in $Q=[0,T]\times [0,1]$.
\end{theorem}

\begin{proof}
We first prove the existence of weak solution for small
time: there exists a small $\delta\in (0,T]$ such that the Cauchy
problem \rref{eq0} and \rref{eq0IC} has a weak solution $\rho\in
C^0([0,\delta];L^1(0,1))$. The idea is to find first the
characteristic curve $\xi=\xi(t)$ passing through $(0,0)$, then
construct a solution to the Cauchy problem.

Let
 \be\label{Omega}
 \Omega_{\delta,M}
  :=\Big\{\xi \in C^0([0,\delta])\colon \xi(0)=0,
  \widetilde\lambda(M)\leq \frac {\xi(s)-\xi(t)}{s-t}\leq \overline\lambda(M),
   \forall s,t\in [0,\delta]\, s>t \Big\},
 \ee
where $\widetilde\lambda, \overline\lambda$ are defined by
\rref{lambda-bound} and
 \be\label{defM}
 M:=\|u\|_{L^1(0,T)}+\|\rho_0\|_{L^1(0,1)}.
 \ee
We point out here that the case $d(M)=0$ (by \rref{lambda-bound}, $\lambda$ is a constant in $[0,M]$)
 is trivial. We only prove Theorem \ref{thm-l1} for the case $d(M)>0$.

We define a map $F:\Omega_{\delta,M} \rightarrow
C^0([0,\delta])$, $\xi\mapsto F(\xi)$, as
 \be \label{Fxi2}
  F(\xi)(t):=\int_0^t\lambda(\int_0^s u(\sigma)
  d\sigma+\int_0^{1-\xi(s)}\rho_0(x)dx)ds, \, \forall \xi\in \Omega_{\delta,M}, \, \forall t \in [0,\delta].
 \ee
It is obvious that $F$ maps into $\Omega_{\delta,M}$ itself if
\be\label{deltasmall}
0<\delta< T \text{ and } \delta < \frac{1}{\overline\lambda(M)}.
 \ee
Now we prove that, if $\delta$ is small enough, $F$
is a contraction mapping on $\Omega_{\delta,M}$ with respect to the
$C^0$ norm defined by
$$
\|\xi\|_{C^0([0,\delta])}:=\sup_{0\leq t\leq \delta}|\xi(t)|.
$$

\vvv

Let  $\xi_1,\xi_2\in \Omega_{\delta,M}$.
We define $\overline\xi_1 \in C^0([0,\delta])$ and $\overline\xi_2 \in C^0([0,\delta])$ by
$\overline\xi_1(t):=\max\{\xi_1(t),\xi_2(t)\}$ and
$\overline\xi_2(t):=\min\{\xi_1(t),\xi_2(t)\}$. By
\rref{lambda-bound} and changing the order of the integrations (see
Figure \ref{figure1}), we have
 \bea\label{Fxi1-Fxi2}
  \lefteqn{|F(\xi_2)(t)-F(\xi_1)(t)| \leq  d(M)
   \int_0 ^t \Big|\int_{1-\xi_1 (s)}^{1-\xi_2 (s)} \rho_0(x) dx \Big| ds}
  \nonumber \\
  &= & d(M)\int_{1-\overline\xi_1(t )}^{1-\overline\xi_2(t )}
    \rho_0(x)(t-\overline\xi_1^{-1}(1-x))dx
 \nonumber \\
  && +d(M)\int_{1-\overline\xi_2(t)}^1 \rho_0(x)
    (\overline\xi_2^{-1}(1-x)- \overline \xi_1^{-1}(1-x)) dx.
 \nonumber \\
 &\leq & d(M)\int_{1-\overline\xi_1(t )}^{1-\overline\xi_2(t)} \rho_0(x) dx
    \cdot (\overline\xi_2^{-1}(\overline \xi_2(t))
    -\overline\xi_1^{-1}(\overline \xi_2(t)))
    \nonumber \\
 && +d(M)\int_{1-\overline\xi_2(t)}^1 \rho_0(x)
    (\overline\xi_2^{-1}(1-x)- \overline \xi_1^{-1}(1-x)) dx
    \nonumber \\
&\leq & d(M)\int_{1-\overline\xi_1(t )}^1 \rho_0(x) dx \cdot
    \sup_{0\leq y\leq \overline \xi_2(t) } (\overline\xi_2^{-1}(y)-\overline\xi_1^{-1}(y)).
 \eea

Using the definitions of $\overline \xi_1,\overline \xi_2$ and of
$\Omega_{\delta,M}$, we obtain that, for every $y\in [0,\overline
\xi_2(t)]$ (see Figure \ref{figure2*}),
 \bea\label{Fxi1-Fxi2-2}
 \lefteqn{ 0\;\leq \;\overline\xi_2^{-1}(y)- \overline \xi_1^{-1}(y)
 }
 \nonumber \\
   &=& \Big(\overline\xi_2^{-1}(y)-\frac{\overline\xi_1^{-1}(y)+\overline\xi_2^{-1}(y)}{2}\Big)
    +\Big(\frac{\overline \xi_1^{-1}(y)+\overline\xi_2^{-1}(y)} {2}-\overline\xi_1^{-1}(y)\Big)
 \nonumber \\
 &\leq &  \frac 1{\widetilde\lambda(M)}
   \Big(y-\overline\xi_2\big (\frac{\overline\xi_1^{-1}(y)+\overline\xi_2^{-1}(y)}{2}\big) \Big)
   +\frac 1{\widetilde\lambda(M)}
    \Big(\overline\xi_1\big (\frac{\overline\xi_1^{-1}(y)+\overline\xi_2^{-1}(y)}{2} \big)-y\Big)
 \nonumber \\
 &= & \frac 1{\widetilde\lambda(M)}
    \Big(\overline\xi_1\big (\frac{\overline\xi_1^{-1}(y)+\overline\xi_2^{-1}(y)}{2} \big)-
    \overline\xi_2\big (\frac{\overline\xi_1^{-1}(y)+\overline\xi_2^{-1}(y)}{2} \big) \Big)
     \nonumber \\
 & \leq & \frac 1{\widetilde\lambda(M)}   \|\xi_1-\xi_2\|_{C^0([0,\delta])}.
 \eea
Therefore,
 \bea
 |F(\xi_2)(t)-F(\xi_1)(t)| & \leq &
  \frac {d(M)}{\widetilde\lambda(M)} \int_{1-\overline \xi_1(t)}^1
   \rho_0(x)dx \cdot \|\xi_1-\xi_2\|_{C^0([0,\delta])}
  \nonumber \\
 &\leq & \frac {d(M)}{\widetilde\lambda(M)} \int_{1-\delta}^1
   \rho_0(x)dx \cdot \|\xi_1-\xi_2\|_{C^0([0,\delta])}.
 \eea

\begin{figure}[h]
 \begin{minipage}[t]{0.45\linewidth}
  \centering
  \includegraphics[width=\textwidth]{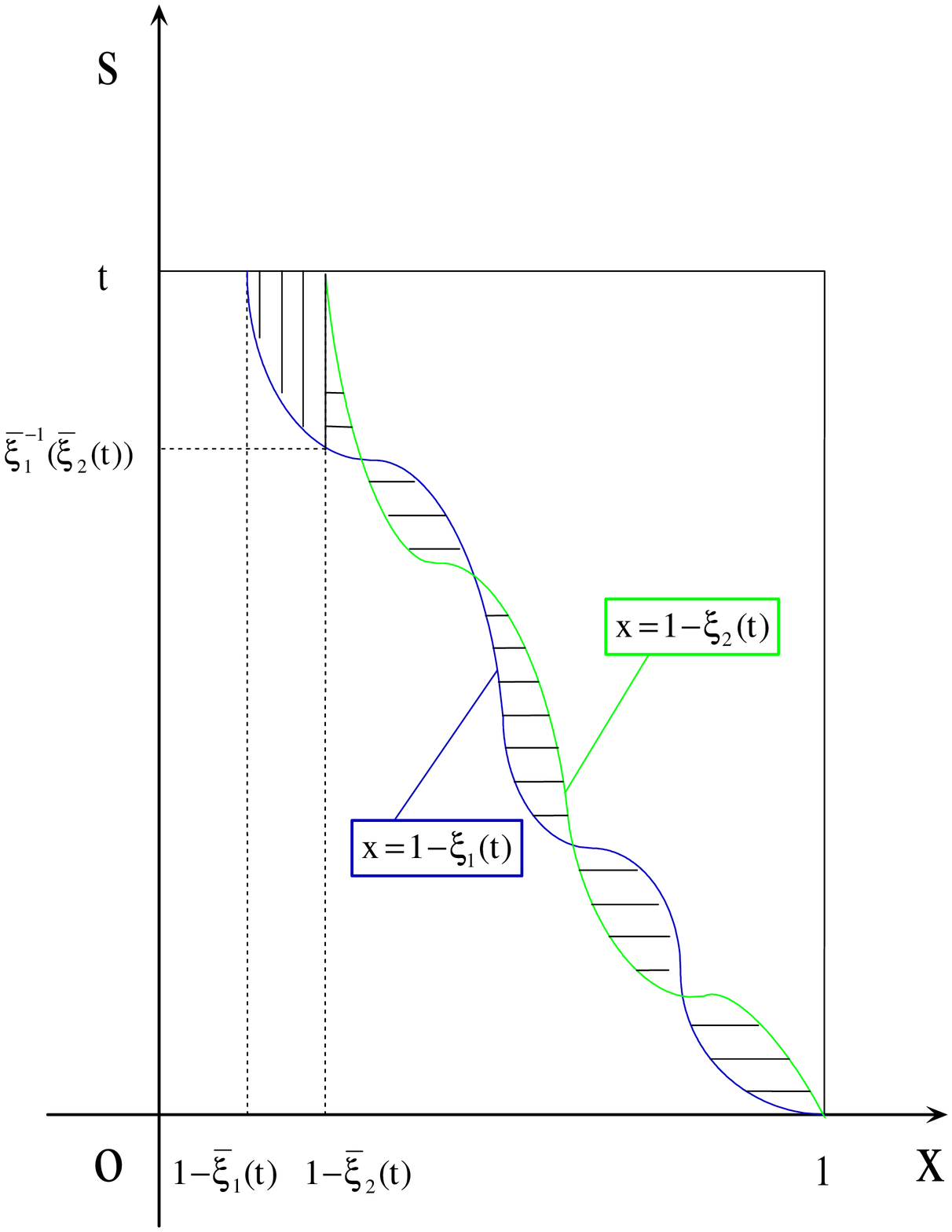}
  \caption{Change order of integrations for $x$ and $s$ in \rref{Fxi1-Fxi2} }
  \label{figure1}
 \end{minipage}
\hspace{0.05\textwidth}
 \begin{minipage}[t]{0.47\linewidth}
  \centering
  \includegraphics[width=\textwidth]{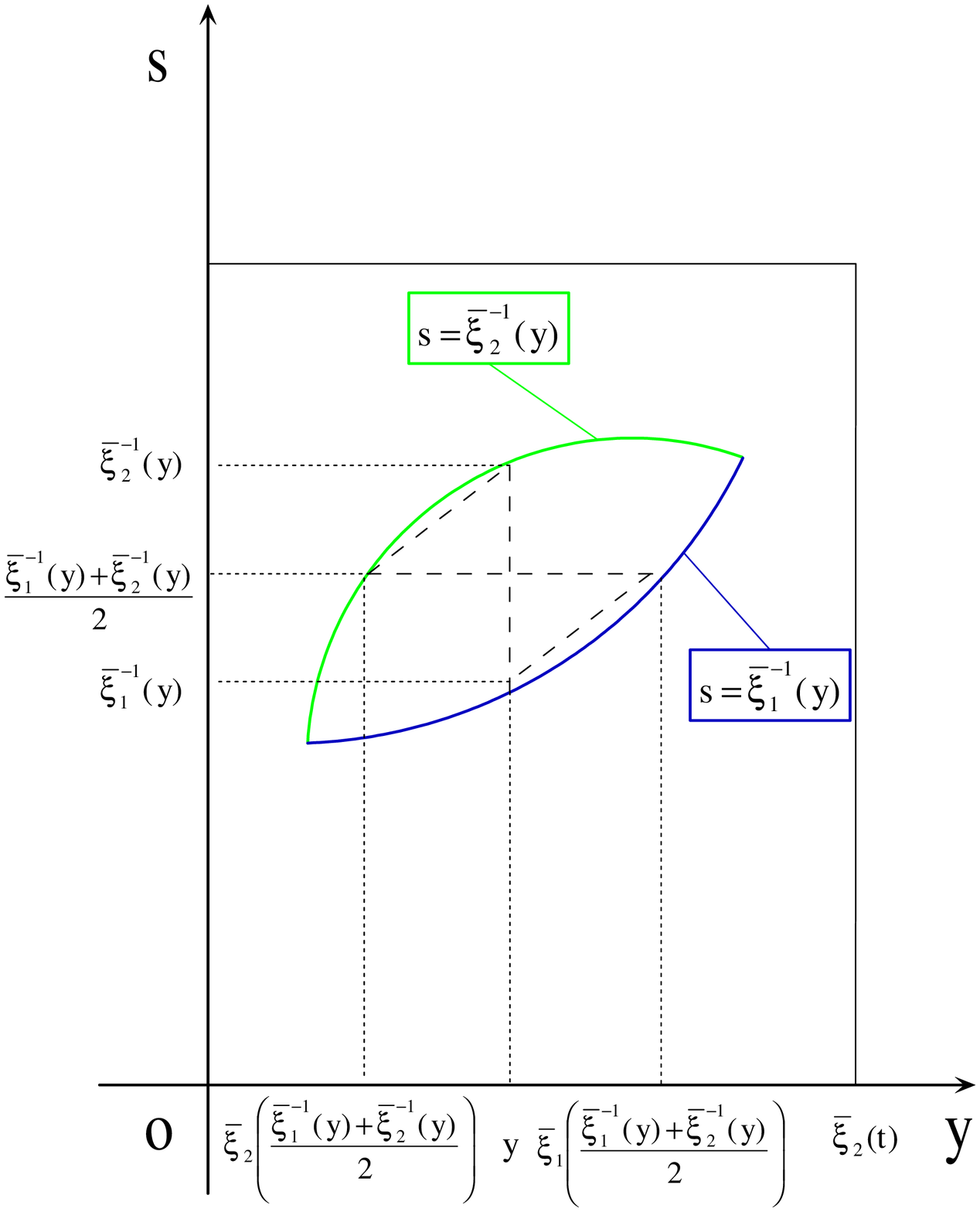}
  \caption{Using property of $\Omega_{\delta,M}$ in \rref{Fxi1-Fxi2-2}}
  \label{figure2*}
\end{minipage}
\end{figure}


 Since $\rho_0\in L^1(0,1)$, we can choose $\delta \in (0,1)$
small enough such that
 \be \label{delta}
  \int_{1-\delta}^1 \rho_0(x)\,dx < \frac {\widetilde\lambda(M)}{2d(M)}.
 \ee
Then
 \be \|F(\xi_1)-F(\xi_2)\|_{C^0([0,\delta])}\leq
 \frac 12 \|\xi_1-\xi_2\|_{C^0([0,\delta])}.
 \ee
By means of the contraction mapping principle, there exists a unique
fixed point $\xi=F(\xi)$ in $\Omega_{\delta,M}$. By \rref{Fxi2}, the
fix point $\xi$ is an increasing function in
$C^1([0,\delta])$, and one has
 \be \label{xi'}
 \xi'(t)=\lambda(\int_0^t u(\sigma) d\sigma+\int_0^{1-\xi(t)}\rho_0(x) dx),
 \quad \forall t\in [0,\delta].
 \ee

Then we define a function $\rho$ by
 \be \label{rho}
 \rho(t,x)=
   \begin{cases}\rho_0(x-\xi(t)),\quad
   & 0\leq \xi(t)\leq x\leq 1,0\leq t\leq \delta,
   \\
   \displaystyle \frac{u(\xi^{-1}(\xi(t)-x))}{\xi'(\xi^{-1}(\xi(t)-x))},\quad
   & 0\leq x\leq \xi(t)\leq 1,0\leq t\leq \delta,
  \end{cases}
 \ee
which is obviously nonnegative almost everywhere. Direct computations
give that, for every $t\in [0,\delta]$,
 \bea \label{wt}
  0\leq  W(t)
 &=&\int_0^1 \rho(t,x) dx
   \nonumber \\
 &=&\int_0^{\xi(t)}\frac{u(\xi^{-1}(\xi(t)-x))}{\xi'(\xi^{-1}(\xi(t)-x))}dx
    +\int_{\xi(t)}^1\rho_0(x-\xi(t))dx
   \nonumber \\
 & =&\int_0^t u(\sigma) d\sigma+\int_0^{1-\xi(t)}\rho_0(y) dy
   \nonumber \\
 & \leq & \|u\|_{L^1(0,T)}+\|\rho_0\|_{L^1(0,1)}= M.
 \eea
Using \rref{lambda-bound},  \rref{xi'} and \rref{wt}, we
obtain the following estimates  of $\xi'$ from above and below:
 \be\label{xi'-bound}
  0<\widetilde\lambda(M)\leq \xi'(t)=\lambda(W(t)) \leq \overline\lambda(M),
  \quad \forall t\in [0,\delta].
 \ee

We now prove that $\rho\in C^0([0,\delta];L^1(0,1))$. For every $s,t\in
[0,\delta]$ with $s\geq t$,
 \bea\label{rhos-rhot}
   \lefteqn{\int_0^1|\rho(s,x)-\rho(t,x)|dx}
   \nonumber\\
 & \leq & \int_0^{\xi(t)}\Big| \frac{u(\xi^{-1}(\xi(s)-x))}{\xi'(\xi^{-1}(\xi(s)-x))}
    -\frac{u(\xi^{-1}(\xi(t)-x))}{\xi'(\xi^{-1}(\xi(t)-x))}\Big|dx
  \nonumber\\
  &&  + \int_{\xi(t)}^{\xi(s)} |\rho(s,x)-\rho(t,x)| dx
    + \int_{\xi(s)}^1 |\rho_0(x-\xi(s))-\rho_0(x-\xi(t))| dx.
 \eea

As for the first term on the right hand side of \rref{rhos-rhot}, we
choose $\{u^n\}_{n=1}^{\infty} \subset C^1([0,T])$ which converges
to $u$ in $L^1(0,T)$, then we have
 \bea \label{rhos-rhot-term-1}
 \lefteqn{ \int_0^{\xi(t)}\Big| \frac{u(\xi^{-1}(\xi(s)-x))}{\xi'(\xi^{-1}(\xi(s)-x))}
    -\frac{u(\xi^{-1}(\xi(t)-x))}{\xi'(\xi^{-1}(\xi(t)-x))}\Big| dx  }
    \nonumber \\
 &\leq &\int_0^{\xi(t)}\Big| \frac{u(\xi^{-1}(\xi(s)-x))}{\xi'(\xi^{-1}(\xi(s)-x))}
    -\frac{u^n(\xi^{-1}(\xi(s)-x))}{\xi'(\xi^{-1}(\xi(s)-x))}\Big| dx
    \nonumber \\
 &&   +\int_0^{\xi(t)}\Big| \frac{u^n(\xi^{-1}(\xi(s)-x))}{\xi'(\xi^{-1}(\xi(s)-x))}
    -\frac{u^n(\xi^{-1}(\xi(t)-x))}{\xi'(\xi^{-1}(\xi(t)-x))}\Big| dx
    \nonumber \\
 &&   +\int_0^{\xi(t)}\Big| \frac{u^n(\xi^{-1}(\xi(t)-x))}{\xi'(\xi^{-1}(\xi(t)-x))}
    -\frac{u(\xi^{-1}(\xi(t)-x))}{\xi'(\xi^{-1}(\xi(t)-x))}\Big| dx
    \nonumber \\
 &\leq & \Big(\int_{\xi^{-1}(\xi(s)-\xi(t))}^s
   +\int_0^t \Big)|u(\sigma)-u^n(\sigma)|  d\sigma
    \nonumber \\
 &&  +\int_0^{\xi(t)}\Big| \frac{u^n(\xi^{-1}(\xi(s)-x))}{\xi'(\xi^{-1}(\xi(s)-x))}
    -\frac{u^n(\xi^{-1}(\xi(t)-x))}{\xi'(\xi^{-1}(\xi(t)-x))}\Big| dx
    \nonumber \\
&\leq & 2 \int_0^T  |u(\sigma)-u^n(\sigma)|  d\sigma
    \nonumber \\
 &&  +\int_0^{\xi(t)}\Big| \frac{u^n(\xi^{-1}(\xi(s)-x))}{\xi'(\xi^{-1}(\xi(s)-x))}
    -\frac{u^n(\xi^{-1}(\xi(t)-x))}{\xi'(\xi^{-1}(\xi(t)-x))}\Big| dx.
 \eea

By \rref{xi'-bound},
 \bea \label{rhos-rhot-term-1.1}
 \lefteqn{\int_0^{\xi(t)} \Big| \frac{u^n(\xi^{-1}(\xi(s)-x))}{\xi'(\xi^{-1}(\xi(s)-x))}
     -\frac{u^n(\xi^{-1}(\xi(t)-x))}{\xi'(\xi^{-1}(\xi(t)-x))}\Big| dx }
    \nonumber\\
 &\leq& \int_0^{\xi(t)}  \Big| \frac{u^n(\xi^{-1}(\xi(s)-x))-u^n(\xi^{-1}(\xi(t)-x))}
    {\xi'(\xi^{-1}(\xi(s)-x))} \Big| dx
   \nonumber\\
 && +\int_0^{\xi(t)} \Big| u^n(\xi^{-1}(\xi(t)-x))\Big(\frac 1{\xi'(\xi^{-1}(\xi(s)-x))}
    -\frac 1{\xi'(\xi^{-1}(\xi(t)-x))} \Big) \Big| dx
   \nonumber\\
 &\leq& C_n |\xi(s)-\xi(t)|
    + C_n \int_0^{\xi(t)} \int_{\xi^{-1}(\xi(t)-x)}^{\xi^{-1}(\xi(s)-x)}
     u(\sigma) d\sigma dx
    \nonumber\\
 &&+ C_n  \int_0^{\xi(t)} \int_{1-\xi(s)+x}^{1-\xi(t)+x}
     \rho_0(y) dy dx,
 \eea
where $C_n$ is a constant independent of $s$ and $t$ but depending
on $u^n$. By changing the order of integrations, we obtain furthermore (see Figure
\ref{figure3})
\bea\label{rhos-rhot-term-1.2}
 \lefteqn{ \int_0^{\xi(t)} \int_{\xi^{-1}(\xi(t)-x)}^{\xi^{-1}(\xi(s)-x)}
     u(\sigma) d\sigma dx  }
    \nonumber\\
 &= & \Big( \int_0^{\xi^{-1}(\xi(s)-\xi(t))} \int_{\xi(t)-\xi(\sigma)}^{\xi(t)}
    +\int_{\xi^{-1}(\xi(s)-\xi(t))}^t \int_{\xi(t)-\xi(\sigma)}^{\xi(s)-\xi(\sigma)}
    + \int_t^s \int_0^{\xi(s)-\xi(\sigma)}\Big) u(\sigma) dx d\sigma
    \nonumber  \\
 &\leq  &  \Big( \int_0^{\xi^{-1}(\xi(s)-\xi(t))}
    +\int_{\xi^{-1}(\xi(s)-\xi(t))}^t +\int_t^s \Big)\,
    u(\sigma) d\sigma \cdot |\xi(s)-\xi(t)|
   \nonumber \\
 &\leq  & \|u\|_{L^1(0,T)}\cdot |\xi(s)-\xi(t)|
 \eea
and (see Figure \ref{figure4})
 \bea\label{rhos-rhot-term-1.3}
 \lefteqn{ \int_0^{\xi(t)} \int_{1-\xi(s)+x}^{1-\xi(t)+x}
     \rho_0(y) dy dx  }
    \nonumber\\
 &= & \Big( \int_{1-\xi(s)}^{1-\xi(t)} \int_0^{\xi(s)-1+y}
    +\int_{1-\xi(t)}^{1-\xi(s)+\xi(t)} \int_{\xi(t)-1+y}^{\xi(s)-1+y}
    + \int_{1-\xi(s)+\xi(t)}^1 \int_{\xi(t)-1+y}^{\xi(t)}\Big) \rho_0(y) dx dy
    \nonumber\\
 &\leq &   \Big(\int_{1-\xi(s)}^{1-\xi(t)}
    +\int_{1-\xi(t)}^{1-\xi(s)+\xi(t)} + \int_{1-\xi(s)+\xi(t)}^1 \Big)
     \rho_0(y) dy \cdot |\xi(s)-\xi(t)|
   \nonumber\\
 &\leq &  \|\rho_0\|_{L^1(0,1)}\cdot |\xi(s)-\xi(t)|.
 \eea

\begin{figure}[h]
 \begin{minipage}[t]{0.47\linewidth}
  \centering
  \includegraphics[width=\textwidth]{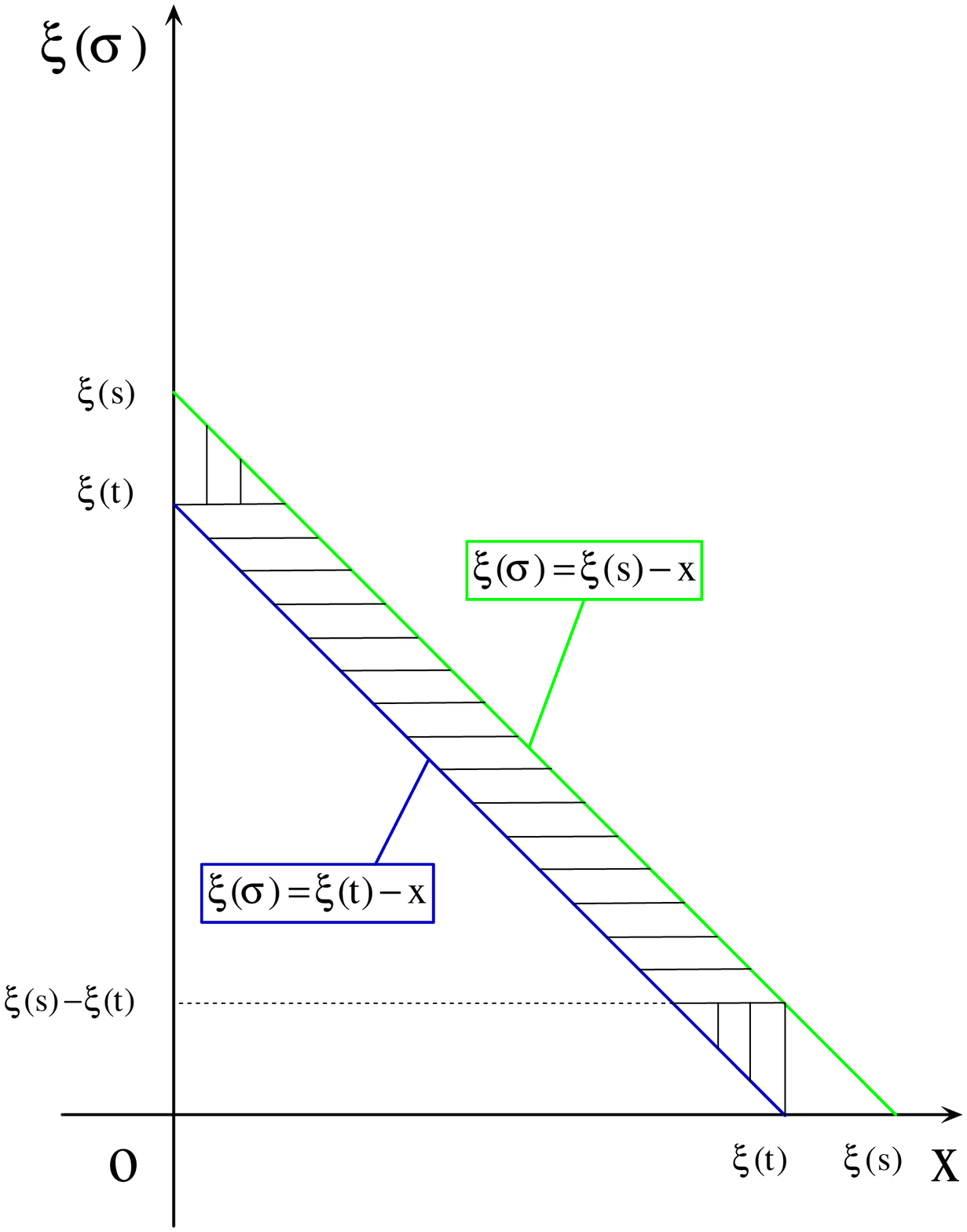}
  \caption{Change order of integration on $\sigma$ and $x$ in
     \rref{rhos-rhot-term-1.2} }
  \label{figure3}
 \end{minipage}
\hspace{0.04\textwidth}
 \begin{minipage}[t]{0.47\linewidth}
  \centering
  \includegraphics[width=\textwidth]{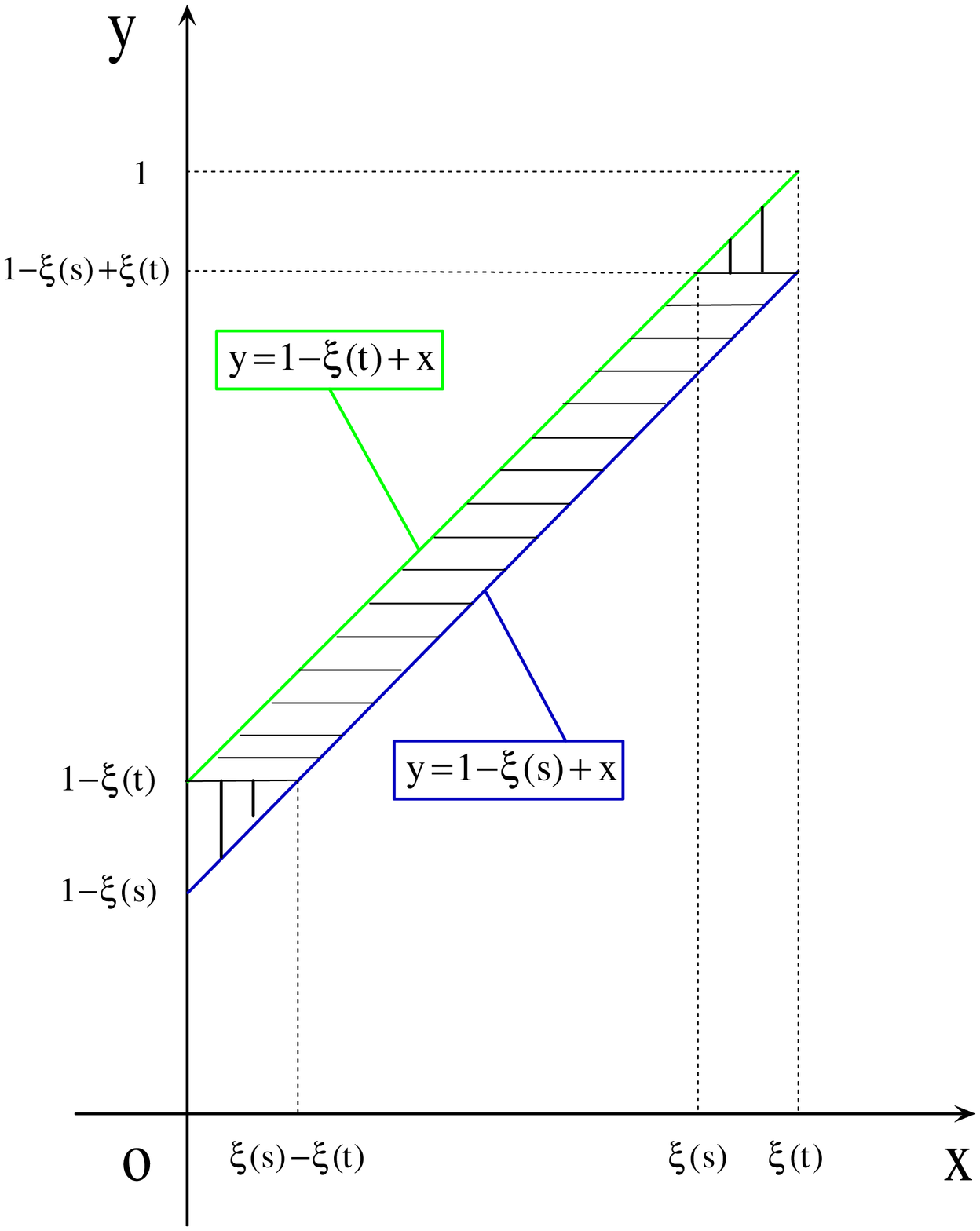}
  \caption{Change order of integration on $y$ and $x$ in
    \rref{rhos-rhot-term-1.3}}
  \label{figure4}
\end{minipage}
\end{figure}


As for the second term on the right hand side of \rref{rhos-rhot}, it is
easy to get that
 \bea\label{rhos-rhot-term-2}
 \lefteqn{ \int_{\xi(t)}^{\xi(s)} |\rho(s,x)-\rho(t,x)| dx
    \leq \int_{\xi(t)}^{\xi(s)} \rho(s,x) dx+ \int_{\xi(t)}^{\xi(s)} \rho(t,x) dx  }
    \nonumber\\
 &= & \int_0^{\xi^{-1}(\xi(s)-\xi(t))} u(\sigma) d\sigma
     +\int_0^{\xi(s)-\xi(t)} \rho_0(y) dy.
 \eea

As for the last term on the right hand side of \rref{rhos-rhot}, we
choose $\{\rho_0^n\}_{n=1}^{\infty} \subset C^1([0,1])$ which
converges to $\rho_0$ in $L^1(0,1)$, then we have
 \bea  \label{rhos-rhot-term-3}
 \lefteqn{\int_{\xi(s)}^1 |\rho_0(x-\xi(s))-\rho_0(x-\xi(t))| dx}
   \nonumber\\
 &\leq & \int_{\xi(s)}^1 |\rho_0(x-\xi(s))-\rho_0^n(x-\xi(s))| dx
    +\int_{\xi(s)}^1 |\rho_0^n(x-\xi(s))-\rho_0^n(x-\xi(t))| dx
    \nonumber\\
 && +\int_{\xi(s)}^1 |\rho_0^n(x-\xi(t))-\rho_0(x-\xi(t))| dx
   \nonumber\\
 &\leq & \Big(\int_0^{1-\xi(s)}+\int_{\xi(s)-\xi(t)}^{1-\xi(t)}\Big)
   |\rho_0(y)-\rho_0^n(y)| dy  + D_n |\xi(s)-\xi(t)|
   \nonumber\\
 &\leq & 2 \int_0^1 |\rho_0(y)-\rho_0^n(y)| dy + D_n  |\xi(s)-\xi(t)|,
 \eea
where $D_n$ is a constant independent of $s$ and $t$ but depending
on $\rho_0^n$.

Using  \rref{defM} together with the estimates \rref{rhos-rhot} to \rref{rhos-rhot-term-3}, we obtain for
any $s,t\in [0,\delta]$ with $s\geq t$,
 \bea\label{rhos-rhot-final}
 \lefteqn{\int_0^1|\rho(s,x)-\rho(t,x)|dx}
    \nonumber\\
 & \leq &   2 \int_0^T  |u(\sigma)-u^n(\sigma)|  d\sigma
    \;\; + \;\; C_n |\xi(s)-\xi(t)|+M |\xi(s)-\xi(t)|
   \nonumber  \\
 &&     +\int_0^{\xi^{-1}(\xi(s)-\xi(t))} u(\sigma) d\sigma
 +\int_0^{\xi(s)-\xi(t)} \rho_0(y) dy
   \nonumber  \\
 &&  + 2 \int_0^1 |\rho_0(y)-\rho_0^n(y)| dy
   + D_n  |\xi(s)-\xi(t)|.
 \eea

We can choose $u^n$ and $\rho_0^n$ such that $\int_0^T
|u(\sigma)-u^n(\sigma)| d\sigma$ and $\int_0^1
|\rho_0(y)-\rho_0^n(y)| dy$ are small as we want. Then according to
\rref{xi'-bound} and the fact that $u\in L^1(0,T)$ and $\rho_0\in
L^1(0,1)$, the right hand side of \rref{rhos-rhot-final} is
sufficiently small if $s$ and $t$ are close enough to each other.
This proves that the function $\rho$ defined by \rref{rho} belongs
to $C^0([0,\delta];L^1(0,1))$.

Next, we  prove that $\rho$ defined by \rref{rho} is a weak solution
to the Cauchy problem \rref{eq0} and \rref{eq0IC}. Let $\delta' \in [0,\delta]$. For any
$\varphi\in C^1([0,\delta']\times [0,1])$ with
$\varphi(\delta',x)\equiv 0$ and $\varphi(t,1)\equiv 0$, let
 \be A:=\int_0^{\delta'} \int_0^1 \rho(t,x)(\varphi_t(t,x)
   +\lambda(W(t))\varphi_x(t,x)) dx dt.
 \ee
Then we have
 \bea
 A &=&\int_0^{\delta'} \int_0^{\xi(t)}
   \frac{u(\xi^{-1}(\xi(t)-x))}{\xi'(\xi^{-1}(\xi(t)-x))}
   (\varphi_t(t,x)+\lambda(W(t))\varphi_x(t,x)) dx dt
 \nonumber\\
 && +\int_0^{\delta'}\int_{\xi(t)}^1 \rho_0(x-\xi(t))
    (\varphi_t(t,x) +\lambda(W(t))\varphi_x(t,x)) dx dt
 \nonumber\\
 &=&\int_0^{\delta'} \int_0^t u(\sigma)(\varphi_t(t,\xi(t)-\xi(\sigma))
     +\lambda(W(t))\varphi_x(t,\xi(t)-\xi(\sigma)) d\sigma dt
  \nonumber\\
 && +\int_0^{\delta'}\int_0^{1-\xi(t)} \rho_0(y)
     (\varphi_t(t,\xi(t)+y) +\lambda(W(t))\varphi_x(t,\xi(t)+y)) dy dt,
\eea
and thus
 \bea
 A &=&\int_0^{\delta'} \int_{\sigma}^{\delta'}
     u(\sigma)\frac{d\varphi(t,\xi(t)-\xi(\sigma))}{dt} dt d\sigma
  \nonumber\\
 && +\Big(\int_0^{1-\xi(\delta')}\int_0^{\delta' }
      +\int_{1-\xi(\delta')}^1\int_0^{\xi^{-1}(1-y)}\Big)
       \rho_0(y)\frac{d\varphi(t,\xi(t)+y)}{dt} dt dy
 \nonumber\\
  &=&-\int_0^{\delta'} u(\sigma)\varphi(\sigma,0)dt
      -\int_0^1\rho_0(y)\varphi(0,y)dy.
 \eea
This proves the existence of weak solutions to the Cauchy problem
\rref{eq0} and \rref{eq0IC} for small time.
\vvv

Now we turn to prove the uniqueness of the weak solution. Let us assume
that $\overline\rho\in C^0([0,\delta];L^1(0,1))$ is a weak solution
to the Cauchy problem \rref{eq0} and \rref{eq0IC}. Then by Lemma \ref{lem-test-function}, for any $\tau\in [0,\delta]$
and $\psi\in C^1([0,\tau]\times[0,1])$ with $\psi(t,1)\equiv 0$,
 \bea \label{test-function}
 \lefteqn{ \int_0^{\tau} \int_0^1\overline\rho(t,x)(\psi_t(t,x)
   +\lambda(\overline W(t))\psi_x(t,x))dxdt +\int_0^{\tau}u(t)\psi(t,0)dt }
  \nonumber\\
 &&-\int_0^1\overline \rho(\tau,x)\psi(\tau,x)dx
    +\int_0^1\rho_0(x)\psi(0,x)dx =0,
 \eea
where $ \overline W(t):=\int_0^1 \overline \rho(t,x)dx $.
\vvv

Let $\overline \xi(t):=\int_0^t \lambda(\overline W(s))ds$ and
$\psi_0\in C_0^1(0,1)$ (i.e. a $C^1$ function with compact support
in $(0,1)$). Then we choose the test function
 \be
 \psi(t,x)=
 \begin{cases}
  \psi_0(\overline \xi(\tau)-\overline \xi(t)+x),\quad
  &0\leq x\leq \overline \xi(t)-\overline \xi(\tau)+1,0\leq t\leq \tau,\\
  0,& 0\leq \overline \xi(t)-\overline \xi(\tau)
  +1 \leq x\leq 1,0\leq  t\leq  \tau,
 \end{cases}
 \ee
which obviously belongs to $C^1([0,\tau] \times [0,1])$ and
satisfies the following backward Cauchy problem:
 \be
 \begin{cases}
 \psi_t+\lambda(\overline W(t))\psi_x=0, \quad
 & 0\leq t\leq \tau, 0\leq x\leq 1,\\
 \psi(\tau,x)=\psi_0(x), &0\leq x\leq 1,\\
 \psi(t,1)=0,& 0\leq t\leq \tau.
 \end{cases}
 \ee

In view of \rref{test-function}, we compute
 \bea
 \lefteqn{\int_0^1 \overline \rho(\tau,x)\psi_0(x) dx}
  \nonumber\\
 &=&\int_0^{\tau} u(t)\psi_0(\overline\xi(\tau)-\overline\xi(t))dt
     +\int_0^{1-\overline\xi(\tau)}\rho_0(x)\psi_0(\overline \xi(\tau)+x)dx
  \nonumber\\
 &=&\int_0^{\overline\xi(\tau)} \frac{u(\overline \xi^{-1}(\overline \xi(\tau)-y))}
     {\overline \xi'(\overline \xi^{-1}(\overline \xi(\tau)-y))}\psi_0(y) dy
     +\int_{\overline \xi(\tau)}^1\rho_0(y-\overline \xi(\tau))\psi_0(y)dy.
 \eea

Since $\psi_0 \in C_0^1(0,1)$ and $\tau \in[0,\delta]$ were arbitrary,
we obtain in $C^0([0,\delta];L^1(0,1))$ that
 \be \label{overrho}
   \overline\rho(t,x)=
 \begin{cases}\rho_0(x-\overline\xi(t)),\quad
   & 0\leq \overline\xi(t)\leq x\leq 1,0\leq t\leq \delta,
   \\\displaystyle
    \frac{u(\overline\xi^{-1}(\overline\xi(t)-x))}
    {\overline\xi'(\overline\xi^{-1}(\overline\xi(t)-x))}, \quad
    & 0\leq x\leq \overline\xi(t)\leq 1,0\leq t\leq \delta,
 \end{cases}
 \ee
which hence gives
 \bea
   \overline  \xi(t)
   &=& \int_0^t\lambda(\int_0^1 \overline\rho(s,x)dx) ds
  \nonumber\\
 &=&\int_0^t\lambda(\int_0^{\overline\xi(s)}
   \frac{u(\overline\xi^{-1}(\overline\xi(t)-x))}
   {\overline\xi'(\overline\xi^{-1}(\overline\xi(t)-x))} dx
   +\int_{\overline\xi(t)}^1\rho_0(x-\overline\xi(s)) dx) ds
 \nonumber\\
 &=&\int_0^t\lambda(\int_0^s u(\sigma)d\sigma
    +\int_0^{1-\overline\xi(s)}\rho_0(y)dy)ds
 \nonumber\\
 &=& F(\overline\xi)(t).
 \eea

It is easy to check that $\overline \xi\in \Omega_{\delta,M}$ when
$\delta$ is small enough, which implies that $\overline \xi=\xi$ since
$\xi$ is the unique fixed point of $F$ in $\Omega_{\delta,M}$ for
$\delta$ small enough, and then $\overline \rho=\rho$ by comparing
\rref{rho} and \rref{overrho}. This gives us the uniqueness of the
weak solution for small time.

Now we suppose that we have solved the Cauchy problem \rref{eq0} and
\rref{eq0IC} to the moment $\tau \in (0,T)$ with the weak solution
$\rho\in C^0([0,\tau];L^1(0,1))$. By \rref{wt}, the following
uniform a priori estimate holds for every $t\in [0,\tau]$:
 \be \label{wt-bound}
 0 \leq W(t)=\int_0^1\rho(t,x)dx \leq  M.
 \ee

Hence we can choose $\delta_0\in (0,T)$ small enough such that \rref{deltasmall} holds and
 $$
 \int_{1-\delta_0}^1\rho(\tau,x)dx
 \leq \frac {\widetilde\lambda(M)}{2d(M)}.
 $$

Applying the previous results on the weak solution for small time,
the weak solution $\rho\in C^0([0,\tau]; L^1(0,1))$ is extended to
the time interval $[\tau,\tau+\delta_0] \cap [\tau,T]$.

Moreover, by \rref{wt-bound} and $\rho_0\in L^1(0,1)$, $u\in
L^1(0,T)$, one can find a suitably small $\delta_0>0$ independent of
$\tau$ such that
 \bea\label{delta0}
 \int_{1-\delta_0}^1\rho(\tau,x)dx
 & \leq & \sup_{t\in [0,T-\frac{\delta_0}{\widetilde\lambda(M)}]}
 \int_t^{t+\frac{\delta_0}{\widetilde\lambda(M)}} u(\sigma) d\sigma
 +\sup_{x\in  [0,1-\delta_0]}
 \int_x^{x+\delta_0} \rho_0(y) dy
 \nonumber \\
 & \leq &  \frac {\widetilde\lambda(M)}{2d(M)}.
 \eea
Step by step, we finally have a unique global weak solution $\rho\in
C^0([0,T];L^1(0,1))$. This concludes the proof of Theorem \ref{thm-l1}.
\end{proof}

\subsection{Remarks}
\begin{remark}\label{remrho}
Let $\rho$ be the weak solution in Theorem \ref{thm-l1}.
 Let $W\in C^0([0,T])$ be defined by $W(t):=\int_0^1\rho(t,x)dx$ and
 let $\xi \in C^1([0,T])$ be defined by requiring
$$
\xi(0)=0, \, \dot \xi (t)=\lambda (W(t)), \, \forall t \in [0,T].
$$
Then, it follows from our proof of Theorem \ref{thm-l1} that
 \be \label{rho-global}
   \rho(t,x)=
 \begin{cases}\rho_0(x-\xi(t)),\quad
    & 0\leq \xi(t)\leq x\leq 1,0\leq t\leq \xi^{-1}(1),
 \\\displaystyle
    \frac{u(\xi^{-1}(\xi(t)-x))}{\xi'(\xi^{-1}(\xi(t)-x))},\quad
    & 0\leq x\leq \xi(t)\leq 1,0\leq t\leq \xi^{-1}(1),
  \\\displaystyle
  \frac{u(\xi^{-1}(\xi(t)-x))}{\xi'(\xi^{-1}(\xi(t)-x))},\quad
  & 0\leq x\leq 1, t\geq \xi^{-1}(1).
 \end{cases}
 \ee
Moreover, $W(t)$ can be expressed as
 \be \label{wt-global}
     W(t)=\int_0^1\rho(t,x)dx=
 \begin{cases}\displaystyle
    \int_0^t u(\sigma) d\sigma
     +\int_0^{1-\xi(t)} \rho_0(y)dy,\quad & 0\leq t\leq \xi^{-1}(1),
   \\\displaystyle
      \int_{\xi^{-1}(\xi(t)-1)}^t u(\sigma) d\sigma
    & t\geq \xi^{-1}(1),
 \end{cases}
 \ee
which implies that
 \be
 0 \leq W(t)=\int_0^1\rho(t,x)dx \leq  M,\quad \forall t\in [0,T]
 \ee
and
 \be\label{xi'-global}
 0<\widetilde\lambda(M)\leq \xi'(t)=\lambda(W(t)) \leq \overline\lambda(M),\quad \forall t\in [0,T].
 \ee
Finally, $W$ is absolutely continuous:
 \be \label{Wt-W't}
 W(t)=W(0)+\int_0^t W'(s)ds
 \ee
with
 \be \label{W't}
     W'(t)=
 \begin{cases}\displaystyle
     u(t)-\xi'(t)\rho_0(1-\xi(t)),\quad
     & 0\leq t\leq \xi^{-1}(1),
  \\ \displaystyle
     u(t)-\frac{\xi'(t)u(\xi^{-1}(\xi(t)-1))}{\xi'(\xi^{-1}(\xi(t)-1))},
    & t\geq \xi^{-1}(1)
\end{cases}
 \ee
and
 \be \label{W't-L1}
 0\leq \int_0^T |W'(t)| dt \leq M.
 \ee
\end{remark}

\begin{remark}
{\bf (Hidden regularity.)} From the definition of the weak solution,
we can expect $\rho\in L^1(0,1;L^1(0,T))$. However, the weak
solution is more regular than expected. In fact, under the assumptions of Theorem \ref{thm-l1},
we have the hidden regularity that $\rho\in C^0([0,1];L^1(0,T))$ so that the function
$t \mapsto \rho(t,x)\in L^1(0,T)$ is well defined for any fixed
$x\in [0,1]$. The proof of the hidden regularity is quite similar to our proof of
 $\rho\in C^0([0,1];L^1(0,T))$ by means of the explicit
expression of $\rho$ (see also  \rref{Wt-W't}-\rref{W't-L1} that we use when $T$ is
large).
\end{remark}

\begin{remark}\label{rem-lp}
If $\rho_0\in L^p(0,1)$ and $u\in L^p(0,T)$ $(p>1)$ are nonnegative
almost everywhere, then the Cauchy problem \rref{eq0} and
\rref{eq0IC} admits a unique weak solution $\rho\in
C^0([0,T];L^p(0,1))\cap C^0([0,1];L^p(0,T))$, which is also
nonnegative almost everywhere in $Q=[0,T]\times [0,1]$. In fact, the
uniqueness of the weak solution comes directly from Theorem
\ref{thm-l1}. And the expression of the solution $\rho\in
C^0([0,T];L^1(0,1))\cap C^0([0,1];L^1(0,T))$ given by
\rref{rho-global} shows that $\rho$ belongs to $C^0([0,T];L^p(0,1))\cap
C^0([0,1];L^p(0,T))$.
\end{remark}

\begin{remark}
If $\rho_0\in C^1([0,1])$ and $u\in C^1([0,T])$ are nonnegative with
 \be
 \begin{cases}
 u(0)-\rho_0(0)=0\\\displaystyle
 u'(0)+\lambda(\int_0^1 \rho_0(x) dx) \rho_0'(0)=0,
 \end{cases}
 \ee
then the Cauchy problem \rref{eq0} and \rref{eq0IC} admits a unique
classical solution $\rho\in C^1([0,T]\times [0,1])$, which is also
nonnegative.
\end{remark}

\section{$L^2$-optimal control for demand tracking problem}

Let $\rho_0\in L^2(0,1)$ be nonnegative almost everywhere and let $T>0$ be given. Let us define
$$L^2_{+}(0,T):=\{u\in L^2(0,T); \text{$u$ is nonnegative almost everywhere}\}.
$$
According to Remark \ref{rem-lp}, for every  $u\in
L^2_{+}(0,T)$, the Cauchy problem \rref{eq0} and \rref{eq0IC} admits a
unique solution $\rho\in C^0([0,T],L^2(0,1)) \cap
C^0([0,1],L^2(0,T))$.

For any fixed {\em demand signal} $y_d\in L^2(0,T)$ and initial data $\rho_0$,
define a functional on $L^2_{+}(0,T)$ by
 \be
 J(u):=\int_0^T |u(t)|^2 dt+ \int_0^T |y(t)-y_d(t)|^2 dt, \, u\in L^2_{+}(0,T),
 \ee
 where
 \be y(t):=\rho(t,1)\lambda(W(t))\ee
is the out-flux corresponding to the in-flux $u\in L^2_{+}(0,T)$ and initial data $\rho_0$.

\begin{theorem}
The infimum of the functional $J$ in $L^2_{+}(0,T)$ is achieved, i.e.,
there exists $u_{\infty}\in L^2_{+}(0,T)$ such that
 \be
 J(u_{\infty})=\inf_{u\in L^2_{+}(0,T)}J(u).
 \ee
\end{theorem}

\begin{proof}
Let $\{u_n\}_{n=1}^{\infty}\subset L^2_{+}(0,T)$ be a minimizing
sequence of the functional $J$, i.e.
 \be
 \lim_{n\rightarrow \infty}J(u_n)=\inf_{u\in L^2_{+}(0,T)}J(u).
 \ee
Then we have
 \be\label{un-yn}
 \|u_n\|_{L^2(0,T)}+\|y_n\|_{L^2(0,T)}\leq C, \quad \forall n\in \mathbb{Z}^+.
 \ee
In \rref{un-yn} and hereafter, we denote by $C$ various constants
which do not depend on $n$.

 The uniform boundedness of $u_n$ in $L^2(0,T)$ shows that there
exists $u_{\infty}\in L^2_{+}(0,T)$ and a subsequence of
$\{u_{n_{k}}\}_{k=1}^{\infty}$ such that $ u_{n_{k}}\rightharpoonup
u_{\infty}$ in $L^2_{+}(0,T)$. For simplicity, we still denote the
subsequence as $\{u_n\}_{n=1}^{\infty}$.

 Let $\rho_n$ be the weak solution to the Cauchy problem of
equation \rref{eq0} with the initial and boundary conditions
 \be
 \begin{cases}
 \rho(t,0)\lambda(W(t))=u_n(t), \quad &0\leq t\leq T,\\
 \rho(0,x)=\rho_0(x),\quad &0\leq x\leq 1.
 \end{cases}
 \ee

Let $W_n:[0,T] \mapsto \mathbb{R}$ and $\xi_n:[0,T] \mapsto
\mathbb{R}$ be defined by
 \be \label{Wn-xin}
 W_n(t):=\int_0^1\rho_n(t,x)dx, \quad
 \xi_n(t):=\int_0^t\lambda(W_n(s))ds.
 \ee

Thus by \rref{wt-global}, we have
 \be \label{xi-n}
\xi_n(t)=\int_0^t\lambda(\int_0^s u_n(\sigma)
d\sigma+\int_0^{1-\xi_n(s)}\rho_0(x)dx)ds,\quad 0\leq t\leq
\min\{\xi_n^{-1}(1),T\}.
 \ee

 In view of \rref{un-yn} and \rref{Wn-xin}, we can derive from \rref{wt-global} that
 \be\label{Wnyn}
 \|W_n\|_{C^0([0,T])}\leq C, \quad \forall n\in \mathbb{Z}^+,
 \ee
which in turn gives with \rref{Wn-xin} that
 \be
 \|\xi_n\|_{C^1([0,T])}\leq C, \quad \forall n\in \mathbb{Z}^+.
 \ee

Moreover, let us point out that $\xi_n'$ is uniformly bounded from
above and below:
 \be
 0 < \widetilde\lambda(\overline C) \leq \xi_n'(t)=\lambda(W_n(t)) \leq
 \overline\lambda(\overline C),\quad \forall t\in [0,T],
 \ \forall n\in \mathbb{Z}^+,
 \ee
where $\widetilde\lambda,\overline \lambda$ are defined by
\rref{lambda-bound} with
 \be
 \overline C:=\sup_{n\in \mathbb{Z}+}
 \|u_n\|_{L^1(0,T)}+\|\rho_0\|_{L^1(0,1)} <\infty.
 \ee

Then it follows from Arzel\`{a}-Ascoli Theorem that there exists
$\overline\xi_{\infty}\in C^0([0,T])$ and a subsequence
$\{\xi_{n_{l}}\}_{l=1}^{\infty}$ such that $ \xi_{n_{l}} \rightarrow
\overline\xi_{\infty}$ in $ C^0([0,T])$. Now we choose the
corresponding subsequence $\{u_{n_{l}}\}_{l=1}^{\infty}$ and again,
denote it as $\{u_n\}_{n=1}^{\infty}$. Thus we have
 \be \label{un-uinfty}
 u_n\rightharpoonup u_{\infty}\quad \text{in}\quad  L^2(0,T),
 \quad \text{as} \ n\rightarrow \infty
 \ee
and
 \be
 \xi_n \rightarrow \overline\xi_{\infty}\quad \text{in}\quad
 C^0([0,T]), \quad \text{as} \ n\rightarrow \infty.
 \ee

Then one has
 \bea
 \lefteqn{\widetilde\lambda(\overline C)|\xi_n^{-1}(x)-\overline\xi_{\infty}^{-1}(x)|
   \leq |\xi_n(\xi_n^{-1}(x))-\xi_n(\overline\xi_{\infty}^{-1}(x))|}
 \nonumber\\
 & &= |x-\xi_n(\overline\xi_{\infty}^{-1}(x))|
    =|\overline\xi_{\infty}(\overline\xi_{\infty}^{-1}(x))
   -\xi_n(\overline\xi_{\infty}^{-1}(x))| \rightarrow 0,\quad
   \text{as}\ n\rightarrow \infty
 \eea
uniformly for $x\in [0,\overline\xi_{\infty}(T))$. Thus we get for
any $x_0\in [0,\overline\xi_{\infty}(T))$,
 \be
 \xi_n^{-1}\rightarrow \overline\xi_{\infty}^{-1} \quad \text{in}\
 C^0([0,x_0]), \quad \text{as} \ n\rightarrow \infty,
 \ee
 and therefore, by passing the limit $n\rightarrow \infty$ in \rref{xi-n},
 \be\label{overxi-infty}
  \overline\xi_{\infty}(t)
  =\int_0^t\lambda(\int_0^s u_{\infty}(\sigma)d\sigma
  +\int_0^{1-\overline\xi_{\infty}(s)}\rho_0(x)dx)ds,\quad
  0\leq t\leq \min\{\overline\xi_{\infty}^{-1}(1),T\}.
 \ee

Let $\rho_{\infty}$ be the weak solution to the Cauchy problem of
equation \rref{eq0} with the initial and boundary conditions
 \be
 \begin{cases}
 \rho(t,0)\lambda(W(t))=u_{\infty}(t), \quad &0\leq t\leq T,\\
 \rho(0,x)=\rho_0(x),\quad &0\leq x\leq 1,
 \end{cases}
 \ee
and denote
 \be\label{WInf-XiInf}
 W_{\infty}(t):=\int_0^1\rho_{\infty}(t,x)dx,\quad
 \xi_{\infty}(t):=\int_0^t\lambda(W_{\infty}(s))ds.
 \ee

We claim that $\xi_{\infty}=\overline \xi_{\infty}$. In fact,
 \be\label{xi-infty}
 \xi_{\infty}(t)
 =\int_0^t\lambda(\int_0^s u_{\infty}(\sigma)d\sigma
 +\int_0^{1-\xi_{\infty}(s)}\rho_0(x)dx)ds,\quad
 0\leq t\leq \min\{\xi_{\infty}^{-1}(1),T\}.
 \ee
As in the proof of Theorem \ref{thm-l1}, there exists $\delta>0$
small enough which is depending only on $u_{\infty}$ and $\rho_0$
such that
 \be
 \xi(t)=F_{\infty}(\xi)(t)
 :=\int_0^t\lambda(\int_0^s u_{\infty}(\sigma) d\sigma
 +\int_0^{1-\xi(s)}\rho_0(x)dx)ds
 \ee
has a unique fixed point in $\Omega_{\delta,\overline C}$ (replacing
$M$ by $\overline C$ in \rref{Omega}). This implies from
\rref{overxi-infty} and \rref{xi-infty} that
$\xi_{\infty}(t)\equiv\overline\xi_{\infty}(t)$ on $ [0,\delta]$.
Moreover, with the help of \rref{delta0}, there exists $\delta_0>0$
independent of $\tau\in (0,T)$ such that if
$\xi_{\infty}(\tau)=\overline\xi_{\infty}(\tau)$ then
$\xi_{\infty}(t)\equiv\overline\xi_{\infty}(t)$ on
$[\tau,\tau+\delta_0]\cap[\tau,T]$.

Therefore
 \be\label{xilim}
 \xi_{\infty} \equiv \overline\xi_{\infty}
 \quad \text{and} \quad
  \xi_n \rightarrow \xi_{\infty}\quad \text{in}\ C^0([0,T]),
  \quad \text{as} \ n\rightarrow \infty,
 \ee
 and  it follows that
 \be\label{xn'lim}
  W_n\rightarrow W_{\infty},
  \quad  \xi_n' \rightarrow \xi_{\infty}'
  \quad \text{in}\ C^0([0,T]),
  \quad \text{as} \ n\rightarrow \infty
  \ee
and, for any $x_0\in[0,\xi_{\infty}(T))$,
 \be \label{xi-1lim}
  \xi_n^{-1}\rightarrow \xi_{\infty}^{-1}
  \quad \text{in}\  C^0([0,x_0]),
  \quad \text{as} \ n\rightarrow \infty.
 \ee

Next prove that $y_n(t)=\lambda(W_n(t))\rho_n(t,1)$ converges
to $y_{\infty}(t)=\lambda(W_{\infty}(t))\rho_{\infty}(t,1)$ weakly
in $L^2(0,T)$. By \rref{un-yn}, $\{y_n\}_{n=1}^{\infty}$ is bounded
in $L^2(0,T)$. Hence, it is suffices to prove that for any $g\in
C^1([0,T])$,
 \be\label{yntoyinf}
 \lim_{n\rightarrow \infty}\int_0^T(y_n(t)-y_{\infty}(t))g(t)dt=0.
 \ee

If $\xi_{\infty}(T)<1$, then $\xi_n(T)<1$ for $n$ large enough. By
\rref{rho-global}, \rref{Wn-xin} and \rref{WInf-XiInf}, for every
$x_0\in [0,\xi_{\infty}(T))$, we have
 \bea\label{Tsmall}
 \lefteqn{\Big|\int_0^T(y_n(t)-y_{\infty}(t))g(t)dt\Big|}
  \nonumber \\
 &=&\Big|\int_0^T(\rho_0(1-\xi_n(t))\lambda(W_n(t))
   -\rho_0(1-\xi_{\infty}(t))\lambda(W_\infty(t))) g(t) dt \Big|
  \nonumber \\
 &=&\Big|\int_{1-\xi_n(T)}^1\rho_0(y)g(\xi_n^{-1}(1-y))dy
    -\int_{1-\xi_{\infty}(T)}^1\rho_0(y)g(\xi_{\infty}^{-1}(1-y))dy\Big|
  \nonumber \\
 &\leq& \Big|\int_{1-\xi_{\infty}(T)}^1\rho_0(y)
    (g(\xi_n^{-1}(1-y))-g(\xi_{\infty}^{-1}(1-y)))dy\Big|
 \nonumber \\
 &+ &     \Big|\int_{1-\xi_n(T)}^{1-\xi_{\infty}(T)}
    \rho_0(y)g(\xi_n^{-1}(1-y)) dy\Big|
 \nonumber \\
 &\leq & \!C\!\! \sup_{0\leq x\leq x_0} |\xi_n^{-1}(x)\!-\!\xi_{\infty}^{-1}(x)|
    + C |\xi_{\infty}(T)\!-\!x_0|^{\frac 12}
    + C |\xi_n(T)\!-\!\xi_{\infty}(T)|^{\frac 12}.
  \eea

By \rref{xilim} and \rref{xi-1lim}, it is easy to get
\rref{yntoyinf} from \rref{Tsmall}.

If $\xi_{\infty}(T)=1$ (i.e., $T=\xi_{\infty}^{-1}(1)$ ), for every
$\tau\in [0,\xi_{\infty}^{-1}(1))$, we have
 \bea\label{T-critical}
 \lefteqn{\Big|\int_0^{\xi_{\infty}^{-1}(1)} (y_n(t)-y_{\infty}(t))g(t)dt\Big|}
  \nonumber\\
 &=& \Big|\int_0^{\tau}(y_n(t)-y_{\infty}(t))g(t)dt
     +\int_{\tau}^{\xi_{\infty}^{-1}(1)} (y_n(t)-y_{\infty}(t))g(t)dt\Big|
  \nonumber\\
 & \leq & \Big|\int_0^{\tau}(y_n(t)-y_{\infty}(t))g(t)dt\Big|
    + C (\xi_{\infty}^{-1}(1)-\tau)^{\frac 12}.
 \eea
Since it is known that for every $\tau\in [0,T)$
$$
\left|\int_0^{\tau}(y_n(t)
-y_{\infty}(t)) g(t) dt \right| \rightarrow 0 \quad \text{ as $n\rightarrow
\infty$},
$$
one has  \rref{yntoyinf} for
$T=\xi_{\infty}^{-1}(1)$ from \rref{T-critical}.

If $\xi_{\infty}(T)>1$, then $\xi_n(T)>1$ for $n$ large enough and  we
have
 \be\label{yn-yInf}
 \int_0^T(y_n(t)-y_{\infty}(t))g(t)dt
 =\Big(\int_0^{\xi_{\infty}^{-1}(1)}
 +\int_{\xi_{\infty}^{-1}(1)}^T\Big) (y_n(t)-y_{\infty}(t))g(t) dt.
 \ee
From the above study, we need only to estimate the last term in
\rref{yn-yInf}. Assuming $\xi_n^{-1}(1)\leq \xi_{\infty}^{-1}(1)$
(the case $\xi_n^{-1}(1)\geq \xi_{\infty}^{-1}(1)$ can be treated
similarly), we get from \rref{rho-global} that
 \bea\label{yn-yInf-2}
 \lefteqn{\Big|\int_{\xi_{\infty}^{-1}(1)}^T(y_n(t)-y_{\infty}(t))g(t)dt\Big|}
  \nonumber\\
 &=&  \Big|\int_{\xi_{\infty}^{-1}(1)}^T(u_n(\xi_n^{-1}(\xi_n(t)-1))
    -u_{\infty}(\xi_{\infty}^{-1}(\xi_{\infty}(t)-1)))g(t)dt\Big|
   \nonumber\\
 &= & \Big|\int_{\xi_n^{-1}(\xi_n(\xi_{\infty}^{-1}(1))-1)}^{\xi_n^{-1}(\xi_n(T)-1)}
    \frac{u_n(\sigma)g(\xi_n^{-1}(\xi_n(\sigma)+1))\xi_n'(\sigma)}
    {\xi_n'(\xi_n^{-1}(\xi_n(\sigma)+1))}d\sigma
   \nonumber\\
 &&-\int_0^{\xi_{\infty}^{-1}(\xi_{\infty}(T)-1)}
     \frac{u_{\infty}(\sigma)g(\xi_{\infty}^{-1}(\xi_{\infty}(\sigma)+1))\xi_{\infty}'(\sigma)}
   {\xi_{\infty}'(\xi_{\infty}^{-1}(\xi_{\infty}(\sigma)+1))}d\sigma \Big|
   \nonumber\\
 &=& \Big|\int_{\tau_n(\xi_{\infty}^{-1}(1))}^{\tau_n(T)}
    \frac{u_n(\sigma)g(\eta_n(\sigma))\xi_n'(\sigma)} {\xi_n'(\eta_n(\sigma))}d\sigma
 \nonumber \\
 &&   -\int_0^{\tau_{\infty}(T)}
     \frac{u_{\infty}(\sigma)g(\eta_{\infty}(\sigma))\xi_{\infty}'(\sigma)}
     {\xi_{\infty}'(\eta_{\infty}(\sigma))}d\sigma \Big|,
 \eea
where we denote
 \bea
  &&\tau_n(t):=\xi_n^{-1}(\xi_n(t)-1),
    \quad \eta_n(t) :=\xi_n^{-1}(\xi_n(t)+1),
   \\
  &&\tau_{\infty}(t):=\xi_{\infty}^{-1}(\xi_{\infty}(t)-1),
    \quad \eta_{\infty}(t):=\xi_{\infty}^{-1}(\xi_{\infty}(t)+1).
 \eea
From \rref{xilim} and \rref{xi-1lim}, we get
 \be\label{taulim}
 \tau_n \rightarrow \tau_{\infty}\quad \text{in} \ C^0([0,T]),
 \quad \text{as} \ n\rightarrow \infty
 \ee
and, for every $t_0\in [0,\tau_{\infty}(T))$,
 \be\label{etalim}
 \eta_n \rightarrow \eta_{\infty}\quad \text{in} \ C^0([t,0_0]),
 \quad \text{as} \ n\rightarrow \infty.
 \ee
Therefore, by \rref{yn-yInf-2}, one has for every $t_0\in
[0,\tau_{\infty}(T))$ and for $n$ large enough,
 \bea
  \lefteqn{\Big|\int_{\xi_{\infty}^{-1}(1)}^T(y_n(t)-y_{\infty}(t))g(t)dt\Big|}
   \nonumber \\
 &\leq & \Big|\Big(\int_{\tau_{\infty}(T)}^{\tau_n(T)}
   -\int_0^{\tau_n(\xi_{\infty}^{-1}(1))}\Big)
   \frac{u_n(\sigma)g(\eta_n(\sigma))\xi_n'(\sigma)}
   {\xi_n'(\eta_n(\sigma))}d\sigma \Big|
   \nonumber \\
 & & +\Big|\int_0^{\tau_{\infty}(T)}
    \Big(\frac{u_n(\sigma)g(\eta_n(\sigma))\xi_n'(\sigma)}
     {\xi_n'(\eta_n(\sigma))}
    -\frac{u_{\infty}(\sigma)g(\eta_{\infty}(\sigma))\xi_{\infty}'(\sigma)}
     {\xi_{\infty}'(\eta_{\infty}(\sigma))}\Big) d\sigma\Big|
   \nonumber \\
  &\leq & C |\tau_n(T)-\tau_{\infty}(T)|^{\frac 12}
      + C |\tau_n(\xi_{\infty}^{-1}(1))|^{\frac 12}
   \nonumber \\
  && + \Big|\int_0^{\tau_{\infty}(T)}
      \frac{u_n(\sigma)\xi_n'(\sigma)}{\xi_n'(\eta_n(\sigma))}
      (g(\eta_n(\sigma))-g(\eta_{\infty}(\sigma)) d\sigma\Big|
   \nonumber \\
  && + \Big|\int_0^{\tau_{\infty}(T)}
       u_n(\sigma)g(\eta_{\infty}(\sigma))
     \Big(\frac{\xi_n'(\sigma)}{\xi_n'(\eta_n(\sigma))}
     -\frac{\xi_{\infty}'(\sigma)}{\xi_{\infty}'(\eta_{\infty}(\sigma))}\Big)
      d\sigma\Big|
  \nonumber \\
  && + \Big|\int_0^{\tau_{\infty}(T)} (u_n(\sigma)-u_{\infty}(\sigma))
      \frac{g(\eta_{\infty}(\sigma))\xi_{\infty}'(\sigma)}
      {\xi_{\infty}'(\eta_{\infty}(\sigma))}d\sigma\Big|
   \nonumber \\
  &\leq & C |\tau_n(T)-\tau_{\infty}(T)|^{\frac 12}
      + C |\tau_n(\xi_{\infty}^{-1}(1))|^{\frac 12}
      + C |\tau_{\infty}(T)-t_0|^{\frac 12}
   \nonumber \\
  &&+ C \sup_{0\leq \sigma \leq t_0}
      |\eta_n(\sigma)-\eta_{\infty}(\sigma)|
      + C \sup_{0\leq \sigma \leq t_0}
      \Big|\frac{\xi_n'(\sigma)}{\xi_n'(\eta_n(\sigma))}
     -\frac{\xi_{\infty}'(\sigma)}
     {\xi_{\infty}'(\eta_{\infty}(\sigma))}\Big|
  \nonumber \\
  && + \Big|\int_0^{\tau_{\infty}(T)} (u_n(\sigma)-u_{\infty}(\sigma))
     \frac{g(\eta_{\infty}(\sigma))\xi_{\infty}'(\sigma)}
    {\xi_{\infty}'(\eta_{\infty}(\sigma))}d\sigma\Big|.
 \eea
By \rref{un-uinfty},\rref{xn'lim}, \rref{taulim}-\rref{etalim} and
the arbitrariness of $t_0 \in[0,\tau_{\infty}(T))$, we have
\rref{yntoyinf} for the case $\xi_{\infty}(T)>1$. This concludes the
proof of \rref{yntoyinf}.

As a result,
 \bea
 J(u_{\infty})&=& \int_0^T|u_{\infty}(t)|^2 dt
    + \int_0^T |y_{\infty}(t)-y_d(t)|^2 dt
  \nonumber\\
 &\leq & \liminf_{n\rightarrow \infty} \int_0^T|u_n(t)|^2 dt
    + \liminf_{n\rightarrow \infty} \int_0^T |y_n(t)-y_d(t)|^2 dt
  \nonumber\\
 & \leq & \liminf_{n\rightarrow \infty} J(u_n)=\lim_{n\rightarrow \infty} J(u_n)
    =\inf_{u\in L^2_{+}(0,T)}J(u)
 \eea
This shows $u_{\infty}$ is a minimizer of $J(u)$ in $L^2_{+}(0,T)$, and
it proves also that $u_n$ tends to $u_{\infty}$ strongly in
$L^2(0,T)$.
\end{proof}

\section{Time-optimal transition between equilibria}
In this section, we focus on the specific model that relates the nonlocal
speed to the total mass according to the assumption \rref{usualspeed}.

It is immediate that constant boundary data $\rho(\cdot,0)=\rho_{\rm in}\geq 0$
eventually drive the state to the equilibrium $\rho \equiv \rho_{\rm in}$.
Together with the symmetry $(t,x,\rho(t,x))\longrightarrow (T-t,1-x,\rho(T-t,1-x))$ of
the conservation law \rref{eq0} this establishes (long-time state) controllability.
Of particular interest is the question of how long it takes to drive the system
from one equilibrium state $\rho_0$ to another
equilibrium state $\rho_1$, compare also the
numerical studies of transfers between equilibria in \cite{lamarca08}.

We first explicitly calculate all quantities for the corresponding
piecewise constant boundary data $\rho(\cdot,0)$, and subsequently
prove that this boundary control is indeed time-optimal.

Suppose $\rho_1\geq \rho_0 \geq 0$ are constant, the initial density is
the equilibrium $\rho(0,x)=\rho_0$ for $x\in (0,1]$, and the
desired terminal density is $\rho(T,x)=\rho_1$ for $x\in [0,1]$
and some minimal $T>0$.
The case $\rho_0\geq \rho_1 \geq 0$ is similar.

A natural choice for the boundary values is $\rho(t,0)=\rho_1$ for $t\geq 0$.
This determines for $0\leq t\leq T$ the control influx and the outflux via
$u(t)=\rho_1\lambda(W(t))$ and $y(t)=\rho_0\lambda(W(t))$,
where $W$ is a solution of the initial value problem
\be
W'(t)={\rho_1-\rho_0\over 1+W(t)},
\rule{4mm}{0mm}
W(0)=\int_0^1 \rho(0,x)\,dx=\rho_0.
\ee
This can be integrated in closed form, yielding
\be
W(t)=-1+\sqrt{(1+\rho_0)^2+2t(\rho_1-\rho_0)}
\ee
and similar expressions for the fluxes and the speed.
All characteristic curves are translations of the solution of the
initial value problem
\be
\xi'(t)=\lambda(W(t))={1\over \sqrt{(1+\rho_0)^2+2t(\rho_1-\rho_0)}},
\rule{4mm}{0mm}
\xi(0)=0.
\ee
which has the explicit solution
\be
\xi(t)={\sqrt{(1+\rho_0)^2+2t(\rho_1-\rho_0)}-(1+\rho_0)
\over \rho_1-\rho_0}.
\ee
The time $T$ to achieve this transition between equilibria
is uniquely determined by $\xi(T)=1$ and evaluates to
\be
\label{exittime1}
T=1+{\rho_0+\rho_1\over 2}.
\ee
In the sequel we prove that this time is indeed minimal.
\vvv

Note that $W$ is a continuous function,  and, in particular $W(T)=\rho_1$.
It is convenient to extend $\rho,u,v$, and $W$ to negative times
by setting $\rho(t,x)=W(t)=\rho_0$ and
$u(t)=y(t)={\rho_0\over 1+\rho_0}=u_0=y_0$ for all $t<0$.
Then $u$ is continuous except for a jump at $t=0$, and
$y$ is continuous except for a jump at $T$. Note that
the height of the jump of $u$ at $t=0$ is larger than
the corresponding jump of $y$ at $T$.
\be
u(0^+)-u(0^-)={\rho_1-\rho_0\over 1+\rho_0}
,\;\;\mbox{ whereas }\;\;
y(T^+)-y(T^-)={\rho_1-\rho_0\over 1+\rho_1}.
\ee
Supposing a jump of the reference demand from
$y_d(t)=\rho_0\lambda(\rho_0)$ for $t<T$ to
$y_d(t)=\rho_1\lambda(\rho_1)$ for $t\geq T$ at this
earliest feasible time, the total backlog
at any $t\geq T$ is, due to the {\em inverse response},
\be
\beta(t)=\int_0^T(y_d(s)-y(s))\,ds=
{\rho_0T\over 1+\rho_0}-\rho_0 \underbrace{\int_0^T \lambda(W(s))\,ds}_{=1}
={\rho_0 T-\rho_0-\rho_0^2\over 1+\rho_0}
\ee
Using the expression \rref{exittime1} for $T$, this simplifies
for $t\geq T$ to
\be
\beta(t)={(\rho_1-\rho_0)\rho_0\over 1+\rho_0}=y_0(\rho_1-\rho_0).
\ee

\begin{figure}[h]
   \centering
   \resizebox{0.5\textwidth}{!}{
   \includegraphics{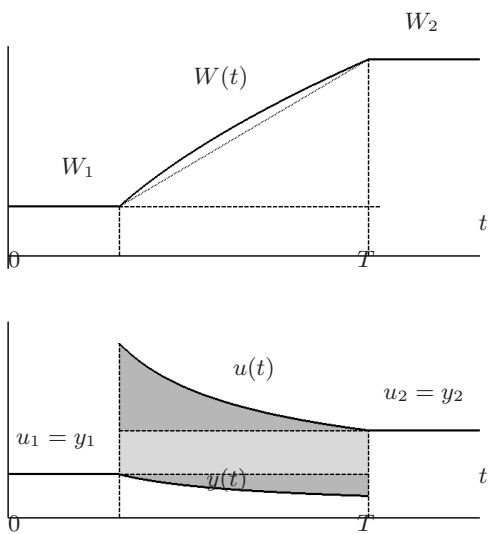}}
   \caption{Total mass, influx, and outflux during optimal transition between equilibria}
   \label{figure5}
\end{figure}
\begin{minipage}{\textwidth}
\begin{picture}(0,0)(0,0)
\put(80,145){{\small $0$}}
\put(212,145){{\small $T$}}
\put(258,159){{\small $t$}}
\put(83,78){{\small $u_1=y_1$}}
\put(165,104){{\small $u(t)$}}
\put(155,62){{\small $y(t)$}}
\put(222,95){{\small $u_2=y_2$}}
\put(100,180){{\small $W_1$}}
\put(150,214){{\small $W(t)$}}
\put(230,235){{\small $W_2$}}
\put(80,45){{\small $0$}}
\put(212,45){{\small $T$}}
\put(258,63){{\small $t$}}
\end{picture}
\end{minipage}

Correspondingly, for $0<t<T$ the total mass $W(t)<\rho_1$
continues to grow, and hence the speed is further decreasing.
Therefore, the influx $u(t)=\rho_1\lambda(W(t)))$ is larger than the
eventual new equilibrium influx $u_1=\rho_1\lambda(\rho_1)$.
The total excess in influx evaluates to
\be
\alpha(T)=\int_0^T (\rho_1\lambda(W(s))-\rho_1\lambda(\rho_1))ds
={(\rho_1-\rho_0)\rho_1\over 1+\rho_1}=u_1(\rho_1-\rho_0)
\ee
Together with the nominal difference $(y_1-y_0)T$ between the accumulated
equilibrium fluxes over the time interval $[0,T]$, these add up the
difference in total mass, compare the three shaded regions in Figure \ref{figure5},
\be
W(T)-W(0)=\int_0^T (u(s)-y(s))\,ds
=\alpha (T)+\beta (T)+(\rho_1\lambda(\rho_1)-\rho_0\lambda (\rho_0))T
=\rho_1-\rho_0.
\ee

While it may seem intuitive that this control is time-optimal, we need to
rigorously prove that it is indeed not possible to improve on this time
by e.g. temporarily increasing the speed via smaller influxes.

\begin{proposition}
\label{mintime}
The minimum time to transfer the state from one equilibrium
$\rho(0,x)=\rho_0$, $x\in (0,1]$ to the equilibrium
$\rho(x,T)=\rho_1>\rho_0$, $x\in [0,1]$ using
influx $u\in L^1([0,\infty),[0,\infty))$
is $T=1+{\rho_0+\rho_1\over 2}$.
\end{proposition}

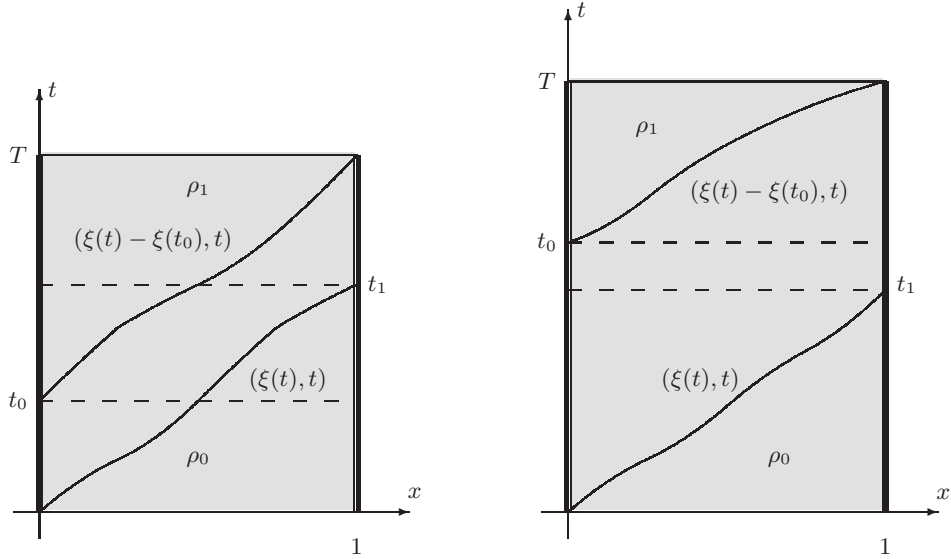
\begin{figure}[htp]
\label{optstep}
\begin{center}
\begin{picture}(400,220)(0,10)
\put(31,21){{\color{litegray}\rule{118pt}{135pt}}}
\put(231,21){{\color{litegray}\rule{118pt}{163pt}}}
\put(20,20){\vector(1,0){150}}
\put(30,10){\vector(0,2){170}}
\put(220,20){\vector(1,0){150}}
\put(230,10){\vector(0,2){200}}
\put(29,20){\line(0,1){135}}
\put(31,20){\line(0,1){135}}
\put(149,20){\line(0,1){135}}
\put(150,20){\line(0,1){135}}
\put(151,20){\line(0,1){135}}
\put(229,20){\line(0,1){163}}
\put(231,20){\line(0,1){163}}
\put(349,20){\line(0,1){163}}
\put(350,20){\line(0,1){163}}
\put(351,20){\line(0,1){163}}
\qbezier(30,20)(45,33)(60,40)
\qbezier(60,40)(75,47)(90,62)
\qbezier(90,62)(105,77)(120,90)
\qbezier(30,62)(45,77)(60,90)
\qbezier(120,90)(135,99)(150,106)
\qbezier(60,90)(75,99)(90,106)
\qbezier(90,106)(105,113)(120,126)
\qbezier(120,126)(135,139)(150,155)
\put(30,155){\line(1,0){120}}
\put(73,123){\makebox(0,0){\small $(\xi(t)-\xi(t_0),t)$}}
\put(124,70){\makebox(0,0){\small $(\xi(t),t)$}}
\put(90,40){\makebox(0,0){\small $\rho_0$}}
\put(90,142){\makebox(0,0){\small $\rho_1$}}
\put(22,155){\makebox(0,0){\small $T$}}
\put(22,62){\makebox(0,0){\small $t_0$}}
\put(158,106){\makebox(0,0){\small $t_1$}}
\put(35,180){\makebox(0,0){\small $t$}}
\put(172,27){\makebox(0,0){\small $x$}}
\put(150,28){\line(0,1){4}}
\put(150,7){\makebox(0,0){\small $1$}}
\multiput(30,62)(12,0){10}{\line(1,0){5}}
\multiput(30,106)(12,0){10}{\line(1,0){5}}
\qbezier(230,20)(245,33)(260,40)
\qbezier(260,40)(275,47)(290,60)
\qbezier(290,60)(305,73)(320,81)
\qbezier(320,81)(335,89)(350,104)
\qbezier(230,122)(245,127)(260,139)
\qbezier(260,139)(275,151)(290,159)
\qbezier(290,159)(305,167)(320,173)
\qbezier(320,173)(335,179)(350,183)
\put(230,183){\line(1,0){120}}
\put(306,141){\makebox(0,0){\small $(\xi(t)-\xi(t_0),t)$}}
\put(280,70){\makebox(0,0){\small $(\xi(t),t)$}}
\put(310,40){\makebox(0,0){\small $\rho_0$}}
\put(260,165){\makebox(0,0){\small $\rho_1$}}
\put(222,183){\makebox(0,0){\small $T$}}
\put(222,122){\makebox(0,0){\small $t_0$}}
\put(358,106){\makebox(0,0){\small $t_1$}}
\put(235,210){\makebox(0,0){\small $t$}}
\put(372,27){\makebox(0,0){\small $x$}}
\put(350,28){\line(0,1){4}}
\put(350,7){\makebox(0,0){\small $1$}}
\multiput(230,104)(12,0){10}{\line(1,0){5}}
\multiput(230,122)(12,0){10}{\line(1,0){5}}
\end{picture}
\caption{Time optimal transfer between equilibrium states}
\end{center}
\end{figure}

\begin{proof}
Suppose $T>0$ and  $\rho(t,0)$ is an integrable function on $[0,T]$
such that the solution of \rref{eq0} satisfies $\rho(T,\cdot)=\rho_1$.

Since $\rho$ is constant along the characteristic curves, there exists
$t_0\in [0,T]$ such that for all $t\in [t_0,T]$, $\rho(t,0)=\rho_1$.
Let $\xi\colon [0,T]\mapsto [0,\infty)$ be the unique function
satisfying $\xi(0)=0$ and $\xi'(t)=\lambda(\int_0^1\rho(t,x)\,dx)$.
Then there exists a unique $t_1\in (0,T]$ such that $\xi(t_1)=1$.

For $0<t<t_0$, $\xi'(t)$ is bounded above by
\be
\xi'(t)={1 \over 1+\int_0^1 \rho(t,x)dx}
\leq{1 \over 1+\int_{\xi(t)}^1 \rho(t,x)dx}
={1 \over 1+\rho_0(1-\xi(t))}.
\ee

Rewrite as $((1+\rho_0)-\rho_0\xi(t))\xi'(t)\leq 1$
and integrate from $t=0$ to $t=t_0$ to obtain a lower
bound for $t_0$.
\be
\label{eqopt1}
(1+\rho_0)\xi(t_0)-{1 \over 2}\rho_0\xi(t_0)^2\leq t_0.
\ee

The primary interest is the case of $t_0 < t_1$.
For $t_0 \leq t\leq t_1$ estimate
\bea
\xi'(t)
& \leq & {1 \over 1+\int_0^{\xi(t)-\xi(t_0)} \rho(t,x)dx+\int_{\xi(t)}^1 \rho(t,x)dx}
\nonumber \\ &=&
{1 \over 1+\rho_1(\xi(t)-\xi(t_0))+\rho_0(1-\xi(t))}.
\eea
and integrate from $t_0$ to $t_1$ to obtain
\be
\label{eqopt2}
(1+\rho_0-\rho_0\xi(t_0))(\xi(t_1)-\xi(t_0))
+{1 \over 2}\,(\rho_1-\rho_0)(\xi(t_1)-\xi(t_0))^2 \leq  t_1-t_0.
\ee

Analogously, for $t_1\leq t\leq T$, the bound
$\xi'(t)\leq 1/(1+\rho_1\,(\xi(t)-\xi(t_0)))$ yields
\be
\label{eqopt3}
\xi(T)-\xi(t_1)+{1 \over 2}\rho_1
\left( (\xi(T)-\xi(t_0))^2 - (\xi(t_1)-\xi(t_0))^2\right)
\leq T-t_1.
\ee

After combining the estimates \rref{eqopt1}, \rref{eqopt2}, and \rref{eqopt3},
elementary simplifications yield
\bea
\label{opt}
T&\geq &
(1+\rho_0)\xi(t_0)-{1 \over 2}\rho_0\xi(t_0)^2
\nonumber \\ &&
+
(1+\rho_0-\rho_0\xi(t_0))(\xi(t_1)-\xi(t_0))
+{1 \over 2}\,(\rho_1-\rho_0)(\xi(t_1)-\xi(t_0))^2
\nonumber \\ &&
+
\xi(T)-\xi(t_1)+{1 \over 2}\rho_1
\left( (\xi(T)-\xi(t_0))^2 - (\xi(t_1)-\xi(t_0))^2\right).
\eea
Noting that $\xi(t_1)=\xi(T)-\xi(t_0)=1$, \rref{opt} simplifies to
\be
T\geq 1+{\rho_0+\rho_1 \over 2} +\xi(t_0).
\ee
This shows that the optimal choice is $t_0=0$, i.e. $\rho(t,0)=\rho_1$ for all $t\geq 0$.
\vvv

It remains to dispose of the case when $t_1 <t_0$.
For $0\leq t \leq t_1$, use $\xi'(t)\leq 1/(1+\rho_0(1-\xi(t)))$
and $\xi(t_1)=1$ to obtain
\be
\label{eqopt4}
t_1 \geq \xi(t_1)(1+\rho_0)-{1 \over 2}\rho_0 \xi(t_1)^2=1+{1 \over 2}\rho_0.
\ee
Similarly, for $t_0 \leq t\leq T$,
use $\xi'(t)\leq 1/(1+\rho_1(\xi(t)-\xi(t_0)))$ and $\xi(T)-\xi(t_0)=1$ to
obtain
\be
\label{eqopt5}
T-t_0\geq (\xi(T)-\xi(t_0))+{1 \over 2}\rho_1 (\xi(T)-\xi(t_0))^2
=1+{1 \over 2}\rho_1.
\ee
Combining  \rref{eqopt4} and \rref{eqopt5} together with
$t_0>t_1$ yields
\be
T=(T-t_0)+(t_0-t_1)+t_1\geq
(1+{1 \over 2}\rho_1)+0+(1+{1 \over 2}\rho_0)
\geq 2+{\rho_0+\rho_1\over 2}.
\ee
This shows that any controls for which $t_1 <t_0$ will perform
even worse than the ones in the first case.
\end{proof}

\section*{Acknowledgements}

The authors thank Fran\c{c}ois Murat for interesting discussions.
The second and third authors thank the Laboratoire Jacques-Louis Lions
at the University Pierre and Marie Curie-Paris VI for its hospitality
and the Foundation Sciences Math\'{e}matiques de Paris for its support.


\begin{thebibliography}{00}

\bibitem{fabio-attset} F. Ancona and A. Marson, On the attainable set
for scalar nonlinear conservation laws with boundary control, SIAM
J. Control Optim., 36 (1998), pp. 290每312 (electronic).

\bibitem{DA-CR-herty5} , Existence theory by front tracking for general
nonlinear hyperbolic systems, Arch. Ration. Mech. Anal., 185 (2007),
pp. 287每340.


\bibitem{DA-CR-hier} D. Armbruster, P. Degond, and C. Ringhofer, A
model for the dynamics of large queuing networks and supply chains,
SIAM J. Appl. Math., 66 (2006), pp. 896每920 (electronic).

\bibitem{fabio-exist} D. Armbruster, D. Marthaler, and C. Ringhofer,
Kinetic and fluid model hierarchies for supply chains, Multiscale
Model. Simul., 2 (2003), pp. 43每61 (electronic).

\bibitem{DA-CR-degond} D. Armbruster, D. Marthaler, C. Ringhofer, K.
Kempf, and T.-C. Jo, A continuum model for a re-entrant factory,
Oper. Res., 54 (2006), pp. 933每950.

\bibitem{DA-CR-thermalized} D. Armbruster and C. Ringhofer, Thermalized
kinetic and fluid models for reentrant supply chains, Multiscale
Model. Simul., 3 (2005), pp. 782每800 (electronic).

\bibitem{piccoli2001} P. Baiti, P. LeFloch, and B. Piccoli, Uniqueness
of classical and nonclassical solutions for nonlinear hyperbolic
systems, J. Differential Equations, 172 (2001), pp. 59每82.

\bibitem{bressanbook} A. Bressan, Hyperbolic systems of conservation
laws, vol. 20 of Oxford Lecture Series in Mathematics and its
Applications, Oxford University Press, Oxford, 2000. The
onedimensional Cauchy problem.

\bibitem{bress-wellposed} A. Bressan, G. Crasta, and B. Piccoli,
Well-posedness of the Cauchy problem for n ℅ n systems of
conservation laws, Mem. Amer. Math. Soc., 146 (2000), pp. viii+134.

\bibitem{shift-diff} A. Bressan and G. Guerra, Shift-differentiability
of the flow generated by a conservation law, Discrete Contin. Dynam.
Systems, 3 (1997), pp. 35每58.

\bibitem{bress97} A. Bressan and P. LeFloch, Uniqueness of weak
solutions to systems of conservation laws, Arch. Rational Mech.
Anal., 140 (1997), pp. 301每317.

\bibitem{piccoli-cars} G. M. Coclite, M. Garavello, and B. Piccoli,
Traffic flow on a road network, SIAM J. Math. Anal., 36 (2005), pp.
1862每1886 (electronic).

\bibitem{herty09} R. Colombo, M. Herty, and M. Mercier, Control of the
continuity equation with a non local flow, preprint, (2009).

\bibitem{coronbook} J.-M. Coron, Control and nonlinearity, vol. 136 of
Mathematical Surveys and Monographs, American Mathematical Society,
Providence, RI, 2007.

\bibitem{coron09} J.-M. Coron, O. Glass, and Z. Wang, Exact boundary
controllability for 1-d quasilinear hyperbolic systems with a
vanishing characteristic speed, preprint, (2009).

\bibitem{herty-network} P. Degond, S. G“ottlich, M. Herty, and A.
Klar, A network model for supply chains with multiple policies,
Multiscale Model. Simul., 6 (2007), pp. 820每837.

\bibitem{herty-supp} M. Herty, A. Klar, and B. Piccoli, Existence of
solutions for supply chain models based on partial differential
equations, SIAM J. Math. Anal., 39 (2007), pp. 160每173.

\bibitem{Horsin} T. Horsin, On the controllability of the Burgers
equation, ESAIM Control Optim. Calc. Var., 3 (1998), pp. 83每95
(electronic).


\bibitem{lamarca08} M. La Marca, D. Armbruster, M. Herty, and C.
Ringhofer, Control of continuum models of production systems,
preprint, (2008).

\bibitem{libook} T. Li, Controllability and Observability for
Quasilinear Hyperbolic Systems, vol. 3 of AIMS Series on Applied
mathematics, 2009.

\bibitem{Li-Rao} T. Li and B. Rao, Exact boundary controllability for
quasi-linear hyperbolic systems, SIAM J. Control Optim., 41 (2003),
pp. 1748每1755 (electronic).

\bibitem{1978-Russell} D. L. Russell, Controllability and
stabilizability theory for linear partial differential equations:
recent progress and open questions, SIAM Rev., 20 (1978), pp.
639每739.

\end{thebibliography}

\end{document}